\documentclass[11pt,leqno]{article}

\usepackage[english]{babel}

\usepackage{amsmath}

\usepackage{amssymb}

\usepackage{xypic}

\usepackage[all]{xy}

\usepackage{latexsym}

\usepackage{amsfonts}

\usepackage{epsf}

\usepackage[dvips]{graphicx}

\usepackage[latin1]{inputenc}

\usepackage[T1]{fontenc}

\numberwithin{equation}{section}

\newcommand{\C}{\mathcal{C}}

\newcommand{\D}{\mathcal{D}}

\newcommand{\T}{\mathcal{T}}

\newcommand{\W}{\mathcal{W}}

\renewcommand{\mod}{\mathrm{Mod}}

\newcommand{\psh}{\mathrm{Psh}}

\newcommand{\ind}{\mathrm{Ind}}

\newcommand{\cov}{\mathrm{Cov}}

\newcommand{\ob}{\mathrm{Ob}}

\newcommand{\op}{\mathrm{Op}}

\newcommand{\rc}{\mathbb{R}\textrm{-}\mathrm{c}}

\newcommand{\CC}{\mathbb{C}}

\newcommand{\R}{\mathbb{R}}

\newcommand{\Z}{\mathbb{Z}}

\newcommand{\N}{\mathbb{N}}

\renewcommand{\P}{\mathcal{P}}

\newcommand{\OO}{\mathcal{O}}

\newcommand{\RR}{\mathcal{R}}

\newcommand{\I}{\mathrm{I}}

\newcommand{\II}{\mathcal{I}}

\newcommand{\osi}{\stackrel{\sim}{\gets}}

\newcommand{\iso}{\stackrel{\sim}{\to}}

\newcommand{\dbtxr}{\mathcal{D}\mathit{b}^t_{X_\mathbb{R}}}

\newcommand{\dbt}{\mathcal{D}\mathit{b}^t}

\newcommand{\db}{\mathcal{D}\mathit{b}}

\newcommand{\ot}{\mathcal{O}^t}

\renewcommand{\th}{\mathit{T}\mathcal{H}\mathit{om}}

\newcommand{\OW}{\OO^\mathrm{w}}
\newcommand{\OWX}{\OO^\mathrm{w}_X}

\newcommand{\CWXR}{\C^{{\infty ,\mathrm{w}}}_{X_\R}}

\newcommand{\CWM}{\C^{{\infty ,\mathrm{w}}}_M}

\newcommand{\wtens}{\overset{\mathrm{w}}{\otimes}}

\newcommand{\rh}{\mathit{R}\mathcal{H}\mathit{om}}

\newcommand{\ho}{\mathcal{H}\mathit{om}}

\newcommand{\Ho}{\mathrm{Hom}}

\newcommand{\Rh}{\mathrm{RHom}}

\newcommand{\id}{\mathrm{id}}

\newcommand{\coker}{\mathrm{coker}}

\renewcommand{\dim}{\textbf{Proof.}}

\newcommand{\qed}{\nopagebreak \phantom{} \hfill $\Box$ \\}

\newcommand{\supp}{\mathrm{supp}}

\newcommand{\imin}[1]{#1^{-1}}

\newcommand{\lind}[1]{\underset{#1}{\underrightarrow{\lim}}}

\newcommand{\Lind}{\underrightarrow{\lim}}  

\newcommand{\indl}[1]{\underset{#1}{``\underrightarrow{\mathrm{lim}}\mbox{''}}}

\newcommand{\Lpro}{\underleftarrow{\lim}}

\newcommand{\lpro}[1]{\underset{#1}{\underleftarrow{\lim}}}

\newcommand{\exs}[3]{0 \to {#1} \to {#2} \to {#3} \to 0}

\newcommand{\lexs}[3]{0 \to {#1} \to {#2} \to {#3}}

\newcommand{\rexs}[3]{{#1} \to {#2} \to {#3} \to 0}

\newtheorem{teo}{Theorem}[subsection]

\newtheorem{df}[teo]{Definition}

\newtheorem{cor}[teo]{Corollary}

\newtheorem{oss}[teo]{Remark}

\newtheorem{prop}[teo]{Proposition}

\newtheorem{lem}[teo]{Lemma}

\newtheorem{nt}[teo]{Notations}

\author{Luca Prelli}
\title{\bf{SHEAVES ON SUBANALYTIC SITES}}
\date{}

\usepackage{latexsym}
\begin{document}


\maketitle

\tableofcontents
\begin{abstract}
In \cite{KS01} the authors introduced the notion of ind-sheaf, and
defined the six Grothendieck operations in this framework. They
defined subanalytic sheaves and they
obtained the formalism of the six Gro\-thendieck operations by
including subanalytic sheaves into the category of ind-sheaves.
The aim of this paper is to give a direct construction of the six
Grothendieck operations
in the framework of subanalytic sites avoiding the heavy theory of
ind-sheaves.
As an application we show how to recover the subanalytic sheaves
$\ot$ and $\OW$ of temperate and Whitney holomorphic functions
respectively.
\end{abstract}
\section*{Introduction}

Let $X$ be a real analytic manifold and $k$ a field. Kashiwara and
Schapira in \cite{KS01} defined and studied the category $\I(k_X)$
of ind-sheaves on $X$. They defined the six Grothendieck
operations in this framework. As a byproduct they also studied the
category $\I_{\rc}(k_X)=\ind(\mod^c_{\rc}(k_X))$, where
$\mod^c_{\rc}(k_X)$ is the category of $\R$-constructible sheaves
on $X$ with compact support, and showed the equivalence with the
category $\mod(k_{X_{sa}})$ of sheaves on the subana\-ly\-tic site
associated to $X$. Then they obtained the formalism of the six
Grothendieck operations by including subanalytic sheaves into the
category of ind-sheaves.

 Our aim in this paper is to give a
direct, self-contained and elementary construction of the six
Grothendieck operations on $\mod(k_{X_{sa}})$, without using the
more sophisticated and much more difficult theory of ind-sheaves.
Indeed, contrarily to the category $\I(k_X)$, the category
$\I_{\rc}(k_X)$ is a Grothendieck category.\\

In more details, the contents of this paper are as follows.

In Section \ref{1} we construct the operations in
$\mod(k_{X_{sa}})$. We start recalling the definitions of the
functors $\rho_*$, $\imin \rho$ and $\rho_!$ of \cite{KS01} and
their properties. We recall the internal operations and the
functors of direct and inverse image (which are well defined on
any site) and we study their relations with $\rho_*$, $\imin \rho$
and $\rho_!$.  We also define the functor $f_{!!}$ of proper
direct image, where the notation $f_{!!}$ follows from the fact
that $f_{!!} \circ \rho_* \not\simeq \rho_* \circ f_!$ in general.
We study its properties and the relations with the others
operations. While the functors $\imin f$ and $\otimes$ are exact,
the functors $\ho$, $f_*$ and $f_{!!}$ are left exact, and we
introduce the subcategory of quasi-injective objects which is
injective with respect to these functors.

In Section \ref{2} we consider the derived category of
$\mod(k_{X_{sa}})$. We start by considering the subcategory
$D^b_{\rc}(k_{X_{sa}})$ consisting of bounded complexes with
$\R$-constructible cohomology, and we prove the equivalence of
derived categories $D^b_{\rc}(k_X) \simeq D^b_{\rc}(k_{X_{sa}})$.
Then we study the derived functors of $\ho$, $f_*$ and $f_{!!}$
and we obtain the usual formulas (projection formula, base change
formula, K\"unneth formula, etc.) in the framework of subanalytic
sites. Using the Brown representability theorem we prove the
existence of a right adjoint to the functor $Rf_{!!}$, denoted by
$f^!$. We calculate the functor $f^!$ by decomposing $f$ as the
composite of a closed embedding and a submersion.

In Section \ref{3} we give some examples of subanalytic sheaves.
Let $\RR$ be a sheaf of rings on $X_{sa}$. We start recalling the
definition of sheaves of $\RR$-modules. When the ring is
$\rho_!\D_X$, where $\D_X$ denotes the sheaf of finite order
differential operators on a complex analytic manifold $X$, we show
how to recover the sheaves of $\rho_!\D_X$-modules $\ot_X$ and
$\OWX$ of temperate and Whitney
holomorphic functions of \cite{KS01} respectively.\\

\noindent \textbf{Acknowledgments.} We would like to thank Prof.
Pierre Schapira who encouraged us to develop a theory of
subanalytic sheaves independent of that of ind-sheaves, and for
his many useful remarks.

\section{Sheaves on subanalytic sites}\label{1}

In the following $X$ will be a real analytic manifold and $k$ a
field. Re\-fe\-ren\-ces are made to \cite{KS} and \cite{Ta94} for
an introduction to sheaves on Grothendieck topologies, to
\cite{KS90} for a complete exposition on classical sheaves and
$\R$-constructible sheaves and to \cite{BM88} and \cite{Lo93} for
the theory of subanalytic sets. The results of $\S$\,\ref{review}
have already been proved in \cite{KS01}, for sake of completeness
we repropose here the proofs.

\subsection{The subanalytic site. Notations and
review}\label{review}

We introduce the subanalytic site and we recall some results of
\cite{KS01} for subanalytic sheaves. We use the notations of \cite{KS01}.\\

Denote by $\op_{sa}(X)$ the category of subanalytic subsets of
$X$. One endows $\op_{sa}(X)$ with the following topology: $S
\subset \op_{sa}(X)$ is a covering of $U \in \op_{sa}(X)$ if for
any compact subset $K$ of $X$ there exists a finite subset $S_0$
of $S$ such that $K \cap \bigcup_{V \in S_0}V=K \cap U$. We will
call $X_{sa}$ the subanalytic site, and for $U \in \op(X_{sa})$ we
denote by $U_{X_{sa}}$ the category $\op(X_{sa}) \cap U$ with the
topology
induced by $X_{sa}$.\\

\begin{oss}\label{UXsa} We use the notation $U_{X_{sa}}$ to stress the
difference from $U_{sa}$, the subanalytic site associated to $U$.
For example, let $X=\R^2$ and $U=\R^2\setminus\{0\}$. Let
$V_n=\{x\in \R^2,\ \ |x|>{1\over n}\}$. Then $\{V_n\}_{n \in \N}
\in \cov(U_{sa})$ but $\{V_n\}_{n\in \N} \notin \cov(U_{X_{sa}})$.
\end{oss}

Let $\mod(k_{X_{sa}})$ denote the category of sheaves on $X_{sa}$.
Then $\mod(k_{X_{sa}})$ is a Grothendieck category, i.e. it admits
a generator and small inductive limits, and small filtrant
inductive limits are exact. In particular as a Grothendieck
category, $\mod(k_{X_{sa}})$ has enough injective objects.\\

\begin{oss}\label{limc} Denote by $\op^c_{sa}(X)$ the category of relatively compact subanalytic open subsets of
$X$. One denotes by $X^c_{sa}$ the category $\op^c_{sa}(X)$ with
the topology induced by $X_{sa}$. The forgetful functor gives an
equivalence of categories
$\mod(k_{X_{sa}})\iso\mod(k_{X^c_{sa}})$.\\
\end{oss}

\begin{prop}\label{UlimU} Let $\{F_i\}_{i \in I}$ be a filtrant inductive
system in $\mod(k_{X_{sa}})$ and let $U \in \op^c(X_{sa})$. Then
$$\lind i\Gamma(U;F_i) \iso \Gamma(U;\lind i F_i).$$
\end{prop}
\dim\ \ By Remark \ref{limc} it is enough to prove the assertion
in the category $\mod(k_{X^c_{sa}})$. Denote by $\indl i F_i$ the
presheaf $V \mapsto \lind i \Gamma(V;F_i)$ on $X^c_{sa}$. Let $U
\in \op^c(X_{sa})$ and let $S$ be a finite covering of $U$. Since
$\lind i$ commutes with finite projective limits we obtain the
isomorphism $(\indl i F_i)(S) \iso \lind i F_i(S)$ and $F_i(U)
\iso F_i(S)$ since $F_i \in \mod(k_{X^c_{sa}})$ for each $i$.
Moreover the family of finite coverings of $U$ is cofinal in
$\cov(U)$. Hence $\indl i F_i \iso (\indl i F_i)^+$. Applying once
again the functor $(\cdot)^+$ we get
$$\indl i F_i \simeq (\indl i F_i)^+ \simeq (\indl i F_i)^{++} \simeq \lind i F_i.$$
Hence applying the functor $\Gamma(U;\cdot)$ we obtain
 the isomorphism $\lind i\Gamma(U;F_i) \iso \Gamma(U;\lind
 i F_i)$ for each $U \in \op^c(X_{sa}).$ \qed \\




There is an easy way to construct sheaves on a subanalytic site

\begin{prop}\label{isiuei} Let $F$ be a presheaf on $X_{sa}^c$ and
assume that
\begin{itemize}
\item[(i)] $F(\varnothing)=0,$
\item[(ii)] For any $U,V \in \op^c(X_{sa})$ the sequence $\lexs{F(U\cap V)}{F(U) \oplus F(V)}{F(U \cap
V)}$is exact.
\end{itemize}
Then $F \in \mod(k_{X_{sa}^c}) \simeq \mod(k_{X_{sa}})$.
\end{prop}
\dim\ \ Let $U \in \op^c(X_{sa})$ and let $\{U_j\}_{j=1}^n$ be a
finite covering of $U$. We have to show that the sequence
$$\lexs{F(U)}{\oplus_{1\leq k\leq n}F(U_k)}{\oplus_{1\leq i<j \leq
n}F(U_{ij})},$$ where the second morphism sends $(s_k)_{1\leq
k\leq n}$ to $(t_{ij})_{1\leq i<j \leq n}$ by
$t_{ij}=s_i|_{U_{ij}}-s_j|_{U_{ij}}$. We shall argue by induction
on $n$. For $n=1$ the result is trivial, and $n=2$ is the
hypothesis. Suppose that the assertion is true for $j \leq n-1$
and set $U'=\bigcup_{1\leq k <n}U_k$. By the induction hypothesis
the following commutative diagram is exact
$$
\xymatrix{&& 0 \ar[d] & 0 \ar[d] \\
0 \ar[r] & F(U) \ar[r] & F(U') \cup F(U_n) \ar[d] \ar[r]
& F(U'\cap U_n) \ar[d] \\
&& \bigoplus_{i<n}F(U_i) \oplus F(U_n) \ar[d] \ar[r] &
\bigoplus_{i<n}F(U_{in}) \\
&& \bigoplus_{i<j<n}F(U_{ij}).}
$$
Then the result follows. \qed \\

Let $\mod_{\rc}(k_X)$  be the abelian category of
$\R$-constructible sheaves on $X$, and consider its subcategory
$\mod^{c}_{\rc}(k_X)$ consisting of sheaves whose support is
compact.\\

We denote by $\rho: X \to X_{sa}$ the natural morphism of sites.
We have functors
\begin{equation}\label{rho}
\xymatrix{\mod^c_{\rc}(k_X)\subset\mod_{\rc}(k_X)\subset\mod(k_X)
\ar@ <2pt> [r]^{\mspace{150mu}\rho_*} &
  \mod(k_{X_{sa}}) \ar@ <2pt> [l]^{\mspace{150mu}\imin \rho}. }
\end{equation}
We will still denote by  $\rho_*$ the restriction of $\rho_*$ to
$\mod_{\rc}(k_X)$ and $\mod^c_{\rc}(k_X)$.\\

\begin{oss}\label{isorho} By Proposition \ref{UlimU} for each $F \in \mod(k_X)$ and $V
\in \op^c(X_{sa})$ one has
$$\Gamma(V;\rho_*F) \simeq \lind U \Gamma(V;\rho_* F_U) \simeq
\Gamma(V;\lind U \rho_*F_U),$$ where $U$ ranges through the family
of relatively compact open subanalytic subsets of $X$. This
implies that $\lind U \rho_* F_U \iso \rho_*F$.
\end{oss}

\begin{oss} The functor $\rho_*$ does not commute with filtrant inductive
limits. For example consider the family $\{V_n\}_{n \in \N}$ of
Remark \ref{UXsa}. We have $\rho_*\lind n k_{V_n} \simeq
\rho_*k_{\R^2\setminus\{0\}}$, while for each $U \in
\op^c_{sa}(\R^2)$ with $0\in \partial U$ we have $\Gamma(U;\lind
n\rho_*k_{V_n}) \simeq \lind n\Gamma(U;\rho_*k_{V_n})=0$.
\end{oss}

\begin{prop}\label{rhoUsa} Let $U$ be an open subanalytic subset of $X$ and consider the constant
sheaf $k_{U_{X_{sa}}} \in \mod(k_{X_{sa}})$. We have
$k_{U_{X_{sa}}} \simeq \rho_*k_U$.
\end{prop}
\dim\ \ Let $F$ be the presheaf defined by $F(V)=k$ if $V \subset
U$, $F(V)=0$ otherwise. This is a separated presheaf and
$k_{U_{X_{sa}}}=F^{++}$. Moreover there is an injective arrow
$F(V) \hookrightarrow \rho_*k_U(V)$ for each $V \in \op(X_{sa})$.
Hence $F^{++} \hookrightarrow \rho_*k_U$ since the functor
$(\cdot)^{++}$ is exact. Let $W \in \op(X_{sa})$ be connected. We
have $F(W) \simeq \rho_*k_U(W) \simeq k$ and then $F^{++} \simeq
\rho_*k_U$ since subanalytic open connected subsets of $X$ form a
basis for
the topology of $X_{sa}$. \qed \\

\begin{prop} One has $\imin \rho \circ \rho_* \iso \id$, in
particular the functor $\rho_*$ is fully faithful.
\end{prop}
\dim\ \ Let $F \in \mod(k_X)$. Every $x \in X$ has a fundamental
neighborhood system consisting of open subanalytic subsets. Hence
we have the chain of isomorphisms
$$(\imin \rho \rho_*F)_x \simeq \lind {x \in U} \imin \rho \rho_*F(U) \simeq \lind {x \in
U} \rho_*F(U) \simeq \lind {x \in U} F(U) \simeq F_x,$$ where $U$
ranges through the family of open subanalytic
neighborhoods of $x$. \qed \\

\begin{prop} The restriction of $\rho_*$ to $\mod_{\rc}(k_X)$ is exact.
\end{prop}
\dim\ \ (i) Let us consider an epimorphism
${G}\twoheadrightarrow{F}$ in $\mod_{\rc}^c(k_X)$, we have to
prove that $\psi:{\rho_*G}\to{\rho_*F}$ is an epimorphism. Let $U
\in \op^c(X_{sa})$ and let $s \in \Gamma(U;\rho_*F) \simeq
\Ho_{k_X}(k_U,F)$. Set $G'=G \times_F k_U$. Then $G' \in
\mod_{\rc}^c(k_X)$ and moreover $G'\twoheadrightarrow k_U$. There
exists a finite $\{U_i\}_{i \in I} \subset \op^c(X_{sa})$ with
$U_i$ connected for each $i$ such that
$\oplus_ik_{U_i}\twoheadrightarrow G'$. The composition $k_{U_i}
\to G' \to k_U$ is given by the multiplication by $a_i \in k$. Set
$I_0=\{k_{U_i};\; a_i \neq 0\}$, we may assume $a_i=1$. We get a
diagram
$$
\xymatrix{\oplus_{i \in I_0}k_{U_i} \ar@{->>} [dr] \ar[r] & G' \ar@{->>} [d] \ar[r] & G \ar@{->>} [d] \\
& k_U \ar[r]^s & F.}
$$
The composition $k_{U_i} \to G' \to G$ defines $t_i \in
\Ho_{k_X}(k_{U_i},G)=\Gamma(U_i;\rho_*G)$. Hence for each $s \in
\Gamma(U;\rho_*F)$ there exists a finite covering $\{U_i\}$ of $U$
and $t_i \in \Gamma(U_i;\rho_*G)$ such that $\psi(t_i)=s|_{U_i}$.
This means that $\psi$ is surjective.

(ii) Let $F \in \mod_{\rc}(k_X)$. By Remark \ref{isorho} $\rho_*F
\simeq \lind U \rho_* F_U$, where $U$ ranges through the family
$\op^c(X_{sa})$. The result follows since $\rho_*$ is exact on
$\mod^c_{\rc}(k_X)$ and filtrant $\Lind$ are exact. \qed
\newline

\begin{nt} Since the functor $\rho_*$ is fully faithful and exact
on the category $\mod_{\rc}(k_X)$, we can identify
$\mod_{\rc}(k_X)$ with its image in $\mod(k_{X_{sa}})$. When there
is no risk of confusion we will write $F$ instead of $\rho_*F$,
for $F \in \mod_{\rc}(k_X)$.
\end{nt}

The following theorem gives a fundamental characterization of
subanalytic sheaves and it will be used systematically in the
following Sections.

\begin{teo}\label{eqlambda} (i) Let $G \in \mod^c_{\rc}(k_X)$ and let
$\{F_i\}$ be a filtrant inductive system in $\mod(k_{X_{sa}})$.
Then we have an isomorphism
$$\lind i \Ho_{k_{X_{sa}}}(\rho_*G,F_i) \iso
\Ho_{k_{X_{sa}}}(\rho_*G,\lind i F_i).$$

(ii) Let $F \in \mod(k_{X_{sa}})$. There exists a small filtrant
inductive system $\{F_i\}_{i \in I}$ in $\mod^c_{\rc}(k_X)$ such
that $F \simeq \lind i \rho_*F_i$.
\end{teo}
\dim\ \ (i) There exists an exact sequence $\rexs{G_1}{G_0}{G}$
with $G_1,G_0$ finite direct sums of constant sheaves $k_U$ with
$U \in \op^c(X_{sa})$. Since $\rho_*$ is exact on
$\mod_{\rc}(k_X)$ and commutes with finite sums, by Proposition
\ref{rhoUsa} we are reduced to prove the isomorphism $\lind i
\Gamma(U;F_i) \iso \Gamma(U;\lind i F_i)$. Then the result follows
from Proposition \ref{UlimU}.

(ii) Let $F \in \mod(k_{X_{sa}})$, and define
\begin{eqnarray*}
I_0 & := & \{(U,s);\; U \in \op^c(X_{sa}),\; s \in \Gamma(U;F)\} \\
G_0 & := & \oplus_{(U,s) \in I_0}\rho_*k_U
\end{eqnarray*}
The morphism $\rho_*k_U \to F$, where the section $1 \in
\Gamma(U;k_U)$ is sent to $s \in \Gamma(U;F)$ defines un
epimorphism $\varphi:G_0 \to F$. Replacing $F$ by $\ker\varphi$ we
construct a sheaf $G_1=\oplus_{(V,t) \in I_1}\rho_*k_V$ and an
epimorphism $G_1\twoheadrightarrow \ker\varphi$. Hance we get an
exact sequence $\rexs{G_1}{G_0}{F}$. For $J_0\subset I_0$ set for
short $G_{J_0}=\oplus_{(U,s)\in J_0}\rho_*k_U$ and define
similarly $G_{J_1}$. Set
$$
J=\{(J_1,J_0);\; J_k\subset I_k,\; J_k \text{ is finite and } {\rm
im} \varphi|_{G_{J_1}} \subset G_{J_0}\}.
$$
The category $J$ is filtrant and $F \simeq \lind {(J_1,J_0)\in J}
\coker(G_{J_1}\to G_{J_0})$. \qed \\

Now we will introduce a left adjoint to the functor $\imin \rho$.

\begin{prop}\label{imineta} Let $F \in \mod(k_{X_{sa}})$, and let $U \in \op(X)$. Then
$$\Gamma(U;\imin \rho F) \simeq \lpro {V \subset\subset U, V \in
\op^c(X_{sa})} \Gamma(V;F)$$
\end{prop}
\dim\ \ By Theorem \ref{eqlambda} we may assume $F=\lind i
\rho_{*}F_i$, with $F_i \in \mod^c_{\rc}(k_X))$. Then $\imin \rho
F \simeq \lind i \imin \rho \rho_{*} F_i \simeq \lind i F_i$. We
have the chain of isomorphisms
$$
\begin{array}{ccccc}
\Gamma(U;\imin \rho F)
 & \hspace{-0cm}\simeq & \hspace{-0.1cm}\lpro{V\subset\subset U
}\Gamma(\overline{V};\imin \rho F)
 & \hspace{-0cm}\simeq & \hspace{-0cm}\lpro{V\subset\subset U
}\Gamma(\overline{V}; \lind i F_i)\\
 & \hspace{-0cm}\simeq & \hspace{-0.1cm}\lpro{V\subset\subset U}\lind{i}\Gamma(\overline{V};F_i)
 & \hspace{-0cm}\simeq & \hspace{-0cm}\lpro{V\subset\subset U}\lind i \hspace{0mm}\Gamma(V;F_i)\\
 & \hspace{-0cm}\simeq & \hspace{-0.1cm}\lpro{V\subset\subset U}\lind i \hspace{0mm}\Gamma(V;\rho_{*}F_i)
 & \hspace{-0cm}\simeq & \hspace{-0cm}\lpro{V\subset\subset U}\Gamma(V;F),
 \end{array}
 $$
where $V \in \op^c(X_{sa})$. The third isomorphism follows since
$\overline{V}$ is compact and the last isomorphism follows from
Proposition \ref{UlimU}. \qed \\

\begin{prop}\label{eta!} The functor $\imin \rho$ admits a left adjoint,
denoted by $\rho_!$. It satisfies
\begin{itemize}
\item[(i)] for $F \in \mod(k_X)$ and $U \in \op(X_{sa})$, $\rho_! F$ is the sheaf associated to the presheaf $U
\mapsto \Gamma(\overline{U};F)$, \\
\item[(ii)] for $U \in \op(X)$ one has $\rho_!k_U \simeq \lind{V\subset\subset U, V \in \op^c(X_{sa})
}k_V$.
\end{itemize}
\end{prop}
\dim\ \ Let $\widetilde{F} \in \psh(k_{X_{sa}})$ be the presheaf
$U \mapsto \Gamma(\overline{U};F)$, and let $G \in
\mod(k_{X_{sa}})$. We will construct morphisms
$$\xymatrix{\Ho_{\psh(k_{X_{sa}})}(\widetilde{F},G) \ar@ <2pt>
[r]^{\ \  \ \xi} &
  \Ho_{k_X}(F,\imin \rho G) \ar@ <2pt> [l]^{\ \ \  \vartheta}}.$$

To define $\xi$, let $\varphi:\widetilde{F} \to G$ and $U \in
\op(X)$. Then the morphism $\xi(\varphi)(U):F(U) \to \imin \rho
G(U)$ is defined as follows

$$F(U) \simeq \lpro{V\subset\subset U, V \in \op^c(X_{sa})
}F(\overline{V}) \stackrel{\varphi}{\longrightarrow}
\lpro{V\subset\subset U, V \in \op^c(X_{sa}) } G(V) \simeq \imin
\rho G(U).$$

On the other hand, let $\psi:F \to \imin \rho G$ and $U \in
\op^c(X_{sa})$. Then the morphism
$\vartheta(\psi)(U):\widetilde{F}(U) \to G(U)$ is defined as
follows

$$\widetilde{F}(U) \simeq \lind {U \subset\subset V \in \op^c(X_{sa})} F(V)
\stackrel{\psi}{\longrightarrow} \lind {U \subset\subset V \in
\op^c(X_{sa})} \imin \rho G(V) \to G(U).$$

By construction one can check that the morphism $\xi$ and
$\vartheta$ are inverse to each others. Then (i) follows from the
chain of isomorphisms
$$\Ho_{k_{X_{sa}}}(\widetilde{F}{}^{++},G) \simeq
\Ho_{\psh(k_{X_{sa}})}(\widetilde{F},G) \simeq
\Ho_{k_{X_{sa}}}(F,\imin \rho G).$$

To show (ii), consider the following sequence of isomorphisms
\begin{eqnarray*}
\Ho_{k_{X_{sa}}}(\rho_!k_U,F) & \simeq & \Ho_{k_X}(k_U,\imin \rho F)\\
 & \simeq & \lpro {V \subset\subset U, V \in \op^c(X_{sa})}
 \Ho_{k_{X_{sa}}}(k_V,F)\\
& \simeq & \Ho_{k_{X_{sa}}}(\lind {V \subset\subset U, V \in
\op^c(X_{sa})}k_V,F),
\end{eqnarray*}
where the second isomorphism follows from Proposition
\ref{imineta}. \qed \\

\begin{prop} The functor $\rho_!$ is exact and commutes with
$\Lind$ and $\otimes$.
\end{prop}
\dim\ \ It follows by adjunction that $\rho_!$ is right exact and
commutes with $\Lind$, so let us show that it is also left exact.
With the notations of Proposition \ref{eta!}, let $F \in
\mod(k_X)$, and let $\widetilde{F} \in \psh(k_{X_{sa}})$ be the
presheaf $U \mapsto \Gamma(\overline{U};F)$. Then $\rho_!F \simeq
\widetilde{F}{}^{++}$, and the functors $F \mapsto
\widetilde{F}$ and $G \mapsto G^{++}$ are left exact.\\
Let us show that $\rho_!$ commutes with $\otimes$. Let $F,G \in
\mod(k_X)$, the morphism
$$F(\overline{U}) \otimes  G(\overline{U}) \to (F\otimes
G)(\overline{U})$$ defines a morphism in $\mod(k_{X_{sa}})$
$$\rho_!F \otimes \rho_!G \to \rho_!(F \otimes G)$$
by Proposition \ref{eta!} (i). Since $\rho_!$ commutes with
$\Lind$ we may suppose that $F=k_U$ and $G=k_V$ and the result
follows from Proposition \ref{eta!} (ii). \qed \\

\begin{prop} The functor $\rho_!$ is fully faithful. In particular
 one has $\imin \rho \circ \rho_! \simeq \id$. Moreover, for $F \in
 \mod(k_X)$ and $G \in \mod(k_{X_{sa}})$ one has
$$\imin \rho \ho(\rho_!F,G) \simeq \ho(F,\imin \rho G).$$
\end{prop}
\dim\ \ For $F,G \in \mod(k_X)$ we have by adjunction
$$\Ho_{k_X}(\imin \rho \rho_!F,G) \simeq \Ho_{k_X}(F,\imin \rho \rho_*G)
\simeq \Ho_{k_X}(F,G).$$ This also implies that $\rho_!$ is fully
faithful, in fact
$$\Ho_{k_{X_{sa}}}(\rho_!F,\rho_!G) \simeq \Ho_{k_X}(F,\imin \rho \rho_!G) \simeq
\Ho_{k_X}(F,G).$$ Now let $K,F \in \mod(k_X)$ and $G \in
\mod(k_{X_{sa}})$, we have
\begin{eqnarray*}
\Ho_{k_X}(K,\imin \rho \ho(\rho_! F,G)) & \simeq & \Ho_{k_{X_{sa}}}(\rho_! K,\ho(\rho_! F,G))\\
 & \simeq & \Ho_{k_{X_{sa}}}(\rho_! K \otimes \rho_! F,G)\\
 & \simeq & \Ho_{k_{X_{sa}}}(\rho_! (K \otimes F),G)\\
 & \simeq & \Ho_{k_X}(K \otimes F,\imin \rho G)\\
 & \simeq & \Ho_{k_X}(K,\ho(F,\imin \rho G))
\end{eqnarray*}
and the result follows.  \qed

\subsection{Operations on the subanalytic site.}

Let $X,Y$ be two real analytic manifolds, and let $f:X \to Y$ be a
real analytic map. This defines a morphism of sites $f:X_{sa} \to
Y_{sa}$. We have a diagram
$$
\xymatrix{X \ar[d]^{\rho} \ar[r]^f & Y \ar[d]^\rho \\
X_{sa} \ar[r]^f & Y_{sa}.}
$$
The following functors are always well defined on a site

\begin{eqnarray*}
\ho & : &  \mod(k_{X_{sa}})^{op} \times \mod(k_{X_{sa}}) \to \mod(k_{X_{sa}}),\\
\otimes & : &  \mod(k_{X_{sa}}^{op})
\times \mod(k_{X_{sa}})  \to \mod(k_{X_{sa}}),\\
f_* & : &  \mod(k_{X_{sa}}) \to \mod(k_{Y_{sa}}),\\
\imin f & : & \mod(k_{Y_{sa}}) \to \mod(k_{X_{sa}}).
\end{eqnarray*}
Let us summarize their properties:
\begin{itemize}
\item the functor $\ho$ is left exact and commutes with $\rho_*$,
\item the functor $\otimes$ is exact and commutes with $\Lind$, $\imin \rho$ and
$\rho_!$,
\item the functor $f_*$ is left exact and commutes with $\rho_*$ and $\Lpro$,
\item the functor $\imin f$ is exact and commutes with $\Lind$, $\otimes$ and $\imin
\rho$,
\item $(\imin f,f_*)$ is a pair of adjoint functors.
\end{itemize}

Let $Z$ be a subanalytic locally closed subset of $X$. As in
classical sheaf theory we define
\begin{eqnarray*}
\Gamma_Z: \mod(k_{X_{sa}}) & \to & \mod(k_{X_{sa}})\\
F & \mapsto & \ho(\rho_*k_Z,F) \\
(\cdot)_Z: \mod(k_{X_{sa}}) & \to & \mod(k_{X_{sa}}) \\
F & \mapsto & F \otimes \rho_*k_Z.
\end{eqnarray*}



We have
\begin{itemize}
\item the functor $\Gamma_Z$ is left exact and commutes with $\rho_*$ and $\Lpro$,
\item the functor $(\cdot)_Z$ is exact and commutes with $\Lind$, $\otimes$ and $\imin
\rho$,
\item $((\cdot)_Z,\Gamma_Z)$ is a pair of adjoint functors.
\end{itemize}

\subsection{$\R$-constructible sheaves on subanalytic
sites}\label{modrcsa}

Let us consider the category $\mod_{\rc}(k_X)$.

\begin{prop}\label{tensrc} Let $F,G \in \mod_{\rc}(k_X)$. Then
$\rho_*(F \otimes G) \simeq \rho_*F \otimes \rho_* G$.
\end{prop}
\dim\ \
We may reduce to the case $F=k_U$, $G=k_V$ with $U,V \in
\op_{sa}(X)$. In this case $\rho_*k_{U \cap V} \simeq \rho_*k_U
\otimes \rho_*k_V$ by Proposition \ref{rhoUsa}. \qed \\


\begin{cor}\label{Zro} Let $F \in \mod_{\rc}(k_X)$, and let $Z$ be a subanalytic
locally closed subset of $X$. Then
$\rho_*F_Z \simeq (\rho_*F)_Z$.
\end{cor}

Let $X,Y$ be two real analytic manifolds, and let $f:X \to Y$ be a
real analytic map.
\begin{prop}\label{iminrc} Let $f:X \to Y$ be a real analytic map.
Let $G \in \mod_{\rc}(k_Y)$. Then $\rho_*\imin f G \simeq \imin f
\rho_*G$.
\end{prop}
\dim\ \
Since the functor $\imin f$ is exact, we may reduce to the case
$G=k_V$, with $V \in \op_{sa}(Y)$. In this case we have $\rho_*
\imin f k_V \simeq \rho_*k_{\imin f (V)} \simeq \imin f
\rho_*k_V$, where the last isomorphism follows from Proposition
\ref{rhoUsa}. \qed \\







We apply the above results to calculate the functor $\ho$ in the
category $\mod(k_{X_{sa}})$.

\begin{prop}\label{hosa} Let $F = \lind i F_i$, with $F_i \in \mod(k_{X_{sa}})$ and let $G = \lind j \rho_*G_j$
with $G_j \in \mod_{\rc}(k_X)$. One has
$$\ho(G,F) \simeq \lpro j \lind i \ho(\rho_*G_j,F_i).$$
\end{prop}
\dim\ \  For each $U \in \op^c(X_{sa})$ one has the isomorphisms
\begin{eqnarray*}
\Gamma(U,\ho(G,F)) & \simeq & \Ho_{k_{X_{sa}}}(G_U,F) \\
& \simeq & \lpro j \lind i\Ho_{k_{X_{sa}}}(\rho_*G_{jU},F_i) \\
& \simeq & \lpro j \lind i\Gamma(U;\ho(\rho_*G_j,F_i)) \\
& \simeq & \Gamma(U;\lpro j \lind i\ho(\rho_*G_j,F_i)).
\end{eqnarray*}
In the second isomorphism we used Corollary \ref{Zro}, and the
last isomorphism follows from Proposition \ref{UlimU} and because
$\Gamma(U;\cdot)$ commutes with $\Lpro$. \qed \\

\begin{cor}\label{horc} Let $F = \lind i \rho_*F_i, G = \lind j \rho_*G_j$
with $F_i,G_j \in \mod_{\rc}(k_X)$. One has
$$\ho(G,F) \simeq \lpro j \lind i \rho_*\ho(G_j,F_i).$$
\end{cor}
\dim\ \ It follows from the fact that $\ho$ commutes with $\rho_*$
and from Proposition \ref{hosa}. \qed \\

\begin{cor} Let $F=\lind i \rho_*F_i$, with $F_i \in
\mod^c_{\rc}(X)$ be a sheaf on $X_{sa}$. Let $Z$ be a subanalytic
locally closed subset of $X$. Then $\Gamma_ZF \simeq \lind i
\rho_* \Gamma_ZF_i.$
\end{cor}

\subsection{Proper direct image on $\mod(k_{X_{sa}})$}

In \cite{KS01} the authors defined the functor $f_{!!}$ of proper
direct image using ind-sheaves.
Here we give a direct construction:
\begin{eqnarray*}
f_{!!}:\mod(k_{X_{sa}}) & \to & \mod(k_{Y_{sa}}) \\
F & \mapsto & \lind U f_*F_U \simeq \lind K f_*\Gamma_KF
\end{eqnarray*}
where $U$ ranges trough the family of relatively compact open
subanalytic subsets of $X$ and $K$ ranges trough the family of
subanalytic compact subsets of $X$. One shall be aware that
$\Lind$ is taken in the category $\mod(k_{Y_{sa}})$. Let $V \in
\op^c(Y_{sa})$. Then $\Gamma(V;f_{!!}F)=\lind U \Gamma(\imin f
(V);F_U) \simeq \lind K \Gamma(\imin f (V);\Gamma_KF),$ where $U$
ranges trough the family of relatively compact open subanalytic
subsets of $X$ and $K$ ranges trough the family of subanalytic
compact subsets of $X$.  If $f$ is proper on $\supp (F)$ then $f_*
\simeq f_{!!}$ and in this case $f_{!!} \circ \rho_* \simeq \rho_*
\circ f_!$.

\begin{oss}
Remark that $f_{!!} \circ \rho_* \neq \rho_* \circ f_!$ in
general. Indeed let $V \in \op^c_{sa}(Y)$, then
\begin{eqnarray*}
\Gamma(V;f_{!!}\rho_*F) & = & \lind K\Ho_{k_X}(k_{\imin f(V)}, \Gamma_KF), \\
\Gamma(V;\rho_*f_!F) & = & \lind Z\Ho_{k_X}(k_{\imin f(V)},
\Gamma_ZF),
\end{eqnarray*}
where $Z$ ranges trough the family of closed subsets of $\imin
f(V)$ such that $f|_Z:Z \to V$ is proper. Then
\begin{eqnarray*}
\Gamma(V;f_{!!}\rho_*F) & = & \{s \in \Gamma(\imin f(V);F);\;
\overline{\supp(s)}
\text{ is compact in X}\}, \\
\Gamma(V;\rho_*f_!F) & = & \{s \in \Gamma(\imin f(V);F);\;
f:\supp(s) \to V \text{ is proper}\}.
\end{eqnarray*}
For example, let $f: \R^2 \to \R$ be the projection on the first
coordinate, and let $V=(a,b) \in \op^c_{sa}(\R)$. Suppose that
$\supp(s)=\{(x,y) \in (a,b) \times \R,\ \
y=\frac{1}{(x-a)(b-x)}\}$. Then $f:\supp(s) \to V$ is proper but
$\overline{\supp(s)}$ is not compact.
\end{oss}

\begin{prop} The functor $f_{!!}$ commutes with filtrant $\Lind$. Moreover $\imin
\rho \circ f_{!!} \simeq f_! \circ \imin \rho$.
\end{prop}
\dim\ \ Let us show that $f_{!!}$ commutes with filtrant $\Lind$.
Let $V \in \op^c_{sa}(Y)$ and let $\{F_i\}_i$ be a filtrant
inductive system in $\mod(k_{X_{sa}})$. Then
\begin{eqnarray*}
\lind K \Ho_{k_{X_{sa}}}(k_{\imin f(V)},\Gamma_K\lind i F_i) &
\simeq & \lind K \Ho_{k_{X_{sa}}}(k_{\imin
f(V) \cap K}, \lind i F_i) \\
& \simeq & \lind {i,K} \Ho_{k_{X_{sa}}}(k_{\imin f(V) \cap K}, F_i) \\
& \simeq & \lind {i,K}\Ho_{k_{X_{sa}}}(k_{\imin f(V)},\Gamma_K F_i) \\
& \simeq & \lind i \Ho_{k_{Y_{sa}}}(k_V,f_{!!}F_i) \\
& \simeq &  \Ho_{k_{Y_{sa}}}(k_V,\lind if_{!!}F_i), \\
\end{eqnarray*}
where the second isomorphism follows from the fact that $k_{\imin
f(V) \cap K} \in \mod_{\rc}^c(k_X)$.\\

Let us show $\imin \rho \circ f_{!!} \simeq f_! \circ \imin \rho$.
Let $F=\lind i \rho_*F_i$. Since $f_{!!}$ commutes with $\Lind$
and $F_i$ has compact support for each $i$ we have $f_{!!}F=\lind
i \rho_*f_!F_i$. We have the chain of isomorphisms
\begin{eqnarray*}
f_!\imin \rho \lind i \rho_*F_i & \simeq & f_!\lind i\imin \rho\rho_*F_i \simeq f_!\lind i F_i \simeq  \lind i f_! F_i \\
& \simeq &   \lind i \imin \rho \rho_* f_! F_i \simeq \imin \rho
\lind i  \rho_* f_! F_i.\ \
\end{eqnarray*}
\qed

\begin{prop} The functor $f_*$ commutes with $\imin \rho$.
\end{prop}
\dim\ \ Let $F \in \mod(k_{X_{sa}})$. Then $f_*F \simeq \lpro K
f_*F_K$, where $K$ ranges through the family of subanalytic
compact subsets of $X$. We have the chain of isomorphisms
\begin{eqnarray*}
f_*\imin \rho F  & \simeq & \lpro K f_*(\imin \rho F)_K \simeq
\lpro K f_{!!}(\imin \rho F)_K \simeq \lpro K f_{!!}\imin \rho
F_K\\
& \simeq & \lpro K \imin \rho f_! F_K \simeq \imin \rho \lpro
Kf_!F_K \simeq \imin \rho \lpro Kf_*F_K\simeq \imin \rho f_*F,
\end{eqnarray*}
where the second and the sixth isomorphism follow from the fact
that $f$ is proper on a compact subset of $X$. \qed \\

\begin{cor} The functor $\imin f$ commutes with $\rho_!$.
\end{cor}
\dim\ \ It follows immediately by adjunction. \qed \\

\begin{prop}\label{bcsa} Let $F \in \mod(k_{X_{sa}})$ and $G \in
\mod(k_{Y_{sa}})$. Then
$$f_{!!}F \otimes G \simeq f_{!!}(F \otimes \imin f G).$$
\end{prop}
\dim\ \ Let $F=\lind i \rho_*F_i$, $G=\lind j \rho_*G_j$. The
functors $\otimes$, $f_{!!}$ and $\imin f$ commute with $\lind i$.
Moreover $\supp(F_i \otimes \imin f G_j)$ is compact for each
$i,j$, hence $f$ is proper on it. Then
\begin{eqnarray*}
f_{!!}\lind i \rho_*F_i \otimes \lind j \rho_*G_j & \simeq & \lind
{i,j}\rho_*(f_!F_i\otimes G_j) \\
& \simeq & \lind {i,j}\rho_* (f_!(F_i\otimes \imin f G_j)) \\
& \simeq & f_{!!}(\lind i\rho_*F_i\otimes \imin f\lind i \rho_*
G_j).
\end{eqnarray*}
In the first isomorphism we used Proposition \ref{tensrc} and in
the last one we used Propositions \ref{tensrc} and \ref{iminrc}.
\qed \\

Now let us consider a cartesian square
$$
\xymatrix{X'_{sa} \ar[r]^{f'} \ar[d]^{g'} & Y'_{sa} \ar[d]^g \\
X_{sa} \ar[r]^f &Y_{sa}}
$$
\begin{prop}\label{pfsa} Let $F \in \mod(k_{X_{sa}})$. Then $\imin g f_{!!}F
\simeq f'_{!!}g'{}^{-1}F$.
\end{prop}
\dim\ \ Let $F=\lind i \rho_*F_i$. All the functors in the above
formula commute with $\lind i$. Moreover since $\supp (F_i)$ is
compact, $f'$ is proper on $\supp (g'{}^{-1}F_i)$ for each $i$.
Then
$$\imin g f_{!!}\lind i \rho_*F_i \simeq \lind i \rho_* \imin g
f_!F_i \simeq \lind i \rho_* f'_!g'{}^{-1}F_i \simeq
f'_{!!}g'{}^{-1}\lind i \rho_*F_i,$$ where the first and the last
isomorphisms follow from Proposition \ref{iminrc}. \qed \\


\begin{prop}\label{nomod} Let $G \in \mod_{\rc}(k_Y)$ and let $F
\in \mod(k_{X_{sa}})$. Then the natural morphism
$$f_{!!}\ho(\imin f G,F) \to \ho(G,f_{!!}F)$$
is an isomorphism.
\end{prop}
\dim\ \ Let us construct the morphism. By adjunction we have
$$\imin f G \otimes \ho(\imin G,F) \to F,$$
hence, using the projection formula we get
$$G \otimes f_{!!}\ho(\imin G, F) \simeq f_{!!}(\imin f G \otimes
\ho(\imin G,F)) \to f_{!!}F,$$ then by adjunction we obtain the
desired morphism. Let us show that it is an isomorphism. We have
the chain of isomorphisms
\begin{eqnarray*}
f_{!!}\ho(\imin f G,F) & \simeq & \lind K f_* \Gamma_K\ho(\imin f
G,F) \\
& \simeq & \lind K f_* \ho(\imin f
G,\Gamma_K F) \\
& \simeq & \lind K \ho(
G,f_*\Gamma_K F) \\
& \simeq & \ho(
G,\lind K f_*\Gamma_K F) \\
& \simeq & \ho(G,f_{!!} F),
\end{eqnarray*}
where the fourth isomorphism follows from Proposition \ref{hosa}.
\qed \\

\subsection{Quasi-injective objects}\label{quinj}

Let us introduce a category which is useful in order to find
acyclic objects with respect to the functors defined in the
previous sections.

\begin{df} An object $F \in \mod(k_{X_{sa}})$ is quasi-injective if the
functor $\Ho_{k_{X_{sa}}}(\cdot,F)$ is exact in
$\mod_{\rc}^c(k_X)$ or, equivalently (see Theorem 8.7.2 of
\cite{KS}) if for each $U,V \in \op^c(X_{sa})$ with $V \subset U$
the restriction morphism $\Gamma(U;F) \to \Gamma(V;F)$ is
surjective.
\end{df}

It follows from the definition that injective sheaves belong to
$\mathcal{J}_{X_{sa}}$. This implies that $\mathcal{J}_{X_{sa}}$
is cogenerating. Moreover the category $\mathcal{J}_{X_{sa}}$ is
stable by filtrant $\Lind$ and $\prod$.


\begin{prop} \label{qinjU} Let  $\exs{F'}{F}{F''}$ be an exact sequence in
$\mod(k_{X_{sa}})$ and assume that $F'$ is quasi-injective. Let $U
\in \op^c(X_{sa})$. Then the sequence
$$\exs{\Gamma(U;F')}{\Gamma(U;F)}{\Gamma(U;F'')}$$
is exact.
\end{prop}
\dim\ \ Let $s'' \in \Gamma(U;F'')$, and let $\{V_i\}_{i=1}^n$ be
a finite covering of $U$ such that there exists $s_i \in
\Gamma(V_i;F)$ whose image is $s''|_{V_i}$. For $n \geq 2$ on $V_1
\cap V_2$ $s_1-s_2$ defines a section of $\Gamma(V_1 \cap V_2;F')$
which extends to $s' \in \Gamma(X;F')$. Replace $s_1$ with
$s_1-s'$. We may suppose that $s_1=s_2$ on $V_1 \cap V_2$. Then
there exists $t \in \Gamma(V_1 \cup V_2)$ such that
$t|_{V_i}=s_i$, $i=1,2$. Thus the induction proceeds. \qed \\

\begin{prop} Let $F',F,F'' \in \mod(k_{X_{sa}})$, and consider the
exact sequence
$$
\exs{F'}{F}{F''}.
$$
Suppose that $F',F \in \mathcal{J}_{X_{sa}}$. Then $F'' \in
\mathcal{J}_{X_{sa}}$.
\end{prop}
\dim\ \ Let $U,V \in \op^c(X_{sa})$ with $V \subset U$ and let us
consider the diagram below
$$ \xymatrix{\Gamma(U;F) \ar[d]^\alpha \ar[r] & \Gamma(U;F'') \ar[d]^\gamma \\
\Gamma(V;F) \ar[r]^\beta & \Gamma(V;F'').} $$  The morphism
$\alpha$ is surjective since $F$ is quasi-injective and $\beta$ is
surjective by Proposition \ref{qinjU}. Then $\gamma$ is
surjective. \qed \\

\begin{teo}\label{rcinj} The family of quasi-injective sheaves is injective with respect to
the functor $\Ho_{k_{X_{sa}}}(G,\cdot)$ for each $G \in
\mod_{\rc}(k_X)$.
\end{teo}
\dim\ \ (i) Let  $\exs{F'}{F}{F''}$ be an exact sequence in
$\mod(k_{X_{sa}})$ and assume that $F',F,F''\in
\mathcal{J}_{X_{sa}}$. Let $G \in \mod_{\rc}^c(k_X)$. We have to
show that the sequence
$$\exs{\Ho_{k_{X_{sa}}}(G,F')}{\Ho_{k_{X_{sa}}}(G,F)}{\Ho_{k_{X_{sa}}}(G,F'')}$$
is exact. $G$ has a resolution
$$0 \to \oplus_{i_1 \in I_1} k_{U_{i_1}} \to \ldots \to \oplus_{i_n \in I_n} k_{U_{i_n}}
\stackrel{\varphi}{\to} G \to 0.$$ Where $I_j$ is finite and
$U_{i_j} \in \op^c_{sa}(X)$ for each $i_j \in I_j$, $j \in
\{1,\ldots ,n\}$. Let us argue by induction on the length $n$ of
the resolution.

(a) If $n=1$, then $G$ is isomorphic to a finite sum $\oplus_i
k_{U_i}$, with $U_i \in \op^c_{sa}(X)$, and the result follows
from Proposition \ref{qinjU}.

(b) Let us show $n-1 \Rightarrow n$. The sequence $\exs{\ker
\varphi}{\oplus_{i_n \in I_n} k_{U_{i_n}}}{G}$ is exact. The sheaf
$\ker \varphi$ belongs to $\mod_{\rc}^c(k_X)$ and it has a
resolution of length $n-1$. We set for short $G_1=\ker\varphi$ and
$G_2=\oplus_{i_n \in I_n} k_{U_{i_n}}$. We get the following
diagram where the columns are exact
$$
\xymatrix{ & 0 \ar[d] & 0 \ar[d] & 0 \ar[d] & \\
0 \ar[r] & \Ho_{k_{X_{sa}}}(G,F') \ar[d] \ar[r] &
\Ho_{k_{X_{sa}}}(G,F)
\ar[d] \ar[r] & \Ho_{k_{X_{sa}}}(G,F'') \ar[d] \ar[r] & 0 \\
0 \ar[r] & \Ho_{k_{X_{sa}}}(G_2,F') \ar[d] \ar[r] &
\Ho_{k_{X_{sa}}}(G_2,F)
\ar[d] \ar[r] & \Ho_{k_{X_{sa}}}(G_2,F'') \ar[d] \ar[r] & 0\\
0 \ar[r] & \Ho_{k_{X_{sa}}}(G_1,F') \ar[d] \ar[r] &
\Ho_{k_{X_{sa}}}(G_1,F)
\ar[d] \ar[r] & \Ho_{k_{X_{sa}}}(G_1,F'') \ar[d] \ar[r] & 0\\
& 0  & 0  & 0  & }
$$

The second row is exact by (i) and the third one is exact by the
induction hypothesis. Hence the top row is exact.

(ii) Let $G \in \mod_{\rc}(k_X)$, let  $\exs{F'}{F}{F''}$ be an
exact sequence in $\mod(k_{X_{sa}})$ with $F'\in
\mathcal{J}_{X_{sa}}$. Let $\{V_n\}_{n \in \N} \in \cov(X_{sa})$
such that $V_n \subset \subset V_{n+1}$. By (i), all the sequences
$$\exs{\Ho_{k_{X_{sa}}}(G_{V_n},F')}{\Ho_{k_{X_{sa}}}(G_{V_n},F)}{\Ho_{k_{X_{sa}}}(G_{V_n},F'')}$$
are exact. Moreover since $F' \in \mathcal{J}_{X_{sa}}$ the
morphism ${\Ho_{k_{X_{sa}}}(G_{V_{n+1}},F')} \to
{\Ho_{k_{X_{sa}}}(G_{V_n},F')}$ is surjective for all $n$. Then by
the Mittag-Leffler property (see Proposition 1.12.3 of
\cite{KS90}) the sequence
$$\exs{\lpro n\Ho_{k_{X_{sa}}}(G_{V_n},F')}{\lpro n\Ho_{k_{X_{sa}}}(G_{V_n},F)}{\lpro n\Ho_{k_{X_{sa}}}(G_{V_n},F'')}$$
is exact. Since $\lpro n\Ho_{k_{X_{sa}}}(G_{V_n},\cdot) \simeq
\Ho_{k_{X_{sa}}}(G,\cdot)$ the result follows. \qed \\

\begin{prop} Let $G \in \mod_{\rc}(k_X)$. Then  quasi-injective sheaves are injective with
respect to the functor $\ho(G,\cdot)$.
\end{prop}
\dim\ \ Let $G \in \mod_{\rc}(k_X)$. It is enough to check that
for each $U \in \op(X_{sa})$ and each exact sequence
$\exs{F'}{F}{F''}$ with $F' \in \mathcal{J}_{X_{sa}}$, the
sequence
$$\exs{\Gamma(U;\ho(G,F'))}{\Gamma(U;\ho(G,F))}{\Gamma(U;\ho(G,F''))}$$
is exact. We have $\Gamma(U,\ho(G,\cdot)) \simeq
\Ho_{k_{X_{sa}}}(G_U,\cdot)$, and quasi-injective objects are
injective with respect to the functor
$\Ho_{k_{X_{sa}}}(G_U,\cdot)$ for each $G \in \mod_{\rc}(k_X)$,
and for each $U \in \op(X_{sa})$. \qed \\

\begin{cor} Quasi-injective sheaves are injective with
respect to the functor $\Gamma_Z$ for each locally closed
subanalytic subset $Z$ of $X$.
\end{cor}

\begin{cor} Let $F \in \mod(k_{X_{sa}})$ be quasi-injective. Then
the functor $\ho(\cdot,F)$ is exact on $\mod_{\rc}(k_X)$.
\end{cor}
\dim\ \ Let $F \in \mod(k_{X_{sa}})$ be quasi-injective. There is
an isomorphism of functors $\Gamma(U;\ho(\cdot,F)) \simeq
\Ho_{k_{X_{sa}}}((\cdot)_U,F)$ for each $U \in \op(X_{sa})$. The
functor $\Ho_{k_{X_{sa}}}((\cdot)_U,F)$
is exact on $\mod_{\rc}(k_X)$ and the result follows. \qed \\

\begin{prop}\label{hoinj} Let $F \in \mod(k_{X_{sa}})$. Then $F$ is
quasi-injective if and only if $\ho(G,F)$ is quasi-injective for
each $G \in \mod_{\rc}(k_X)$.
\end{prop}
\dim\ \ (i) Let $F$ be quasi-injective, and let $G \in
\mod_{\rc}(k_X)$.
We have $\Ho_{k_{X_{sa}}}(\cdot,\ho(G,F)) \simeq
\Ho_{k_{X_{sa}}}(\cdot \otimes G,F)$, and $\Ho_{k_{X_{sa}}}(\cdot
\otimes G,F)$ is exact on $\mod^c_{\rc}(k_X)$.

(ii) Suppose that $\ho(G,F)$ is quasi-injective for each $G \in
\mod_{\rc}(k_X)$. The result follows by setting $G=k_X$.
\qed \\

\begin{cor} The functor $\Gamma_Z$ send quasi-injective
objects to quasi-injective objects for each locally closed
subanalytic subset $Z$ of $X$.
\end{cor}

Let $f:X \to Y$ be a morphism of real analytic manifolds.

\begin{prop} Quasi-injective sheaves are injective with
respect to the functor $f_*$. The functor $f_*$  sends
quasi-injective objects to quasi-in\-jec\-ti\-ve objects.
\end{prop}
\dim\ \ (i) Let us consider $V \in \op(Y_{sa})$. There is an
isomorphism of functors $\Gamma(V;f_*(\cdot)) \simeq \Gamma(\imin
f(V);\cdot)$. It follows from Proposition \ref{rcinj} that
$\mathcal{J}_{X_{sa}}$ is injective with respect to the functor
$\Gamma(\imin f(V);\cdot) \simeq \Ho_{k_{X_{sa}}}(k_{\imin f
(V)},\cdot)$ for any $V \in \op(Y_{sa})$.

(ii) Let $F \in \mathcal{J}_{X_{sa}}$.
For each $G \in \mod_{\rc}^c(k_Y)$ we have
$\Ho_{k_{Y_{sa}}}(G,f_*F) \simeq \Ho_{k_{X_{sa}}}(\imin fG,F)$.
Since $\imin f$ is exact and sends $\mod_{\rc}^c(k_Y)$ to
$\mod_{\rc}(k_X)$, Proposition \ref{rcinj} implies that the
functor $\Ho_{k_{X_{sa}}}(\imin f(\cdot),F)$ is e\-xact on
$\mod^c_{\rc}(k_Y)$. \qed \\

\begin{prop} The family of quasi-injective sheaves is
$f_{!!}$-injective. The functor $f_{!!}$  sends quasi-injective
objects to quasi-in\-jec\-ti\-ve objects.
\end{prop}
(i) Let  $\exs{F'}{F}{F''}$ be an exact sequence in
  $\mod(k_{X_{sa}})$ and assume that $F'\in \mathcal{J}_{X_{sa}}$. We have to check that
 the sequence $\exs{f_{!!}F'}{f_{!!}F}{f_{!!}F''}$ is exact.
Since $F'\in \mathcal{J}_{X_{sa}}$, we have $\Gamma_KF'\in
\mathcal{J}_{X_{sa}}$. Moreover $\mathcal{J}_{X_{sa}}$ is
injective with respect to $\Gamma_K$ and $f_*$. This implies that
the sequence
$$\exs{f_*\Gamma_KF'}{f_*\Gamma_KF}{f_*\Gamma_KF''}$$
is exact. Applying the exact functor $\lind K$ we find that the
sequence
$$\exs{\lind K f_*\Gamma_KF'}{\lind Kf_*\Gamma_KF}{\lind Kf_*\Gamma_KF''}$$
is exact.

(ii) Let $K$ be a compact subanalytic subset of $X$. The functors
$\Gamma_K$ and $f_*$ send quasi-injective objects to
quasi-injective objects, then $f_*\Gamma_KF \in
\mathcal{J}_{Y_{sa}}$. Since $\mathcal{J}_{Y_{sa}}$ is stable by
filtrant $\Lind$,  the result follows. \qed \\

Let $S$ be a closed subanalytic subset of $X$ and let $i_S:S
\hookrightarrow X$ be the closed embedding. Let $F = \lind i\rho_*
F_i \in \mod(k_{X_{sa}})$ with $F_i \in \mod^c_{\rc}(k_X)$. We
have $F_S \simeq \lind i \rho_* F_{iS} \simeq \lind i \rho_*
i_{S*} \imin {i_S} F_i \simeq i_{S*} \imin {i_S}F$.

\begin{lem}\label{WVS} Let $S$ be a closed subanalytic subset of $X$ and
let $U \in \op^c(X_{sa})$. Let $F \in \mod(k_{X_{sa}})$. Then
$\Gamma(U;F_S) \simeq \lind {V \supset S \cap U} \Gamma(V;F)$,
with $V \in \op^c(X_{sa})$.
\end{lem}
\dim\ \ Let $F \in \mod(k_{X_{sa}})$. Then $F \simeq \lind i
\rho_*F_i$ with $F_i \in \mod_{\rc}^c(k_X)$. We have the chain of
isomorphisms
\begin{eqnarray*}
\Gamma(U;F_S) & \simeq & \lind i \Gamma(U;F_{iS}) \\
& \osi & \lind {i,V \supset S \cap U} \Gamma(V;F_i) \\
& \simeq & \lind {V \supset S \cap U} \Gamma(V;F),
\end{eqnarray*}
where $V$ ranges through the family of relatively compact open
subanalytic subsets of $X$ containing $S \cap U$. The second
isomorphism follows since the $F_i$ is $\R$-constructible for each
$i$. \qed \\


\begin{prop}\label{FSqinj} Let $S$ be a closed subanalytic subset of $X$ and let
$F \in \mod(k_{X_{sa}})$ be quasi-injective. Then $F_S$ is
quasi-injective.
\end{prop}
\dim\ \ Let $U,V \in \op^c(X_{sa})$ with $V \subset U$. Since $F$
is quasi-injective and inductive limits are right exact, the
morphism $\lind {U' \supset S \cap U}\Gamma(U'; F) \to \lind {V'
\supset S \cap V}\Gamma(V'; F)$ with $V',U' \in \op^c(X_{sa})$, is
surjective. Hence by Lemma \ref{WVS} the morphism $\Gamma(U;F_S)
\to \Gamma(V;F_S)$ is surjective and the
result follows. \qed \\


Recall that $F \in \mod(k_X)$ is c-soft if the natural morphism
$\Gamma(X;F) \to \Gamma(K,F)$ is surjective for each compact $K
\subset X$. If $F$ is c-soft and $Z$ is a locally closed subset of
$X$, then $F_Z$ is c-soft. Moreover c-soft sheaves are
$\Gamma(U;\cdot)$-injective for each $U \in \op(X)$.

\begin{prop}\label{qsoft} Let $F \in \mod(k_{X_{sa}})$ be quasi-injective. then
$\imin \rho F$ is c-soft.
\end{prop}
\dim\ \ Let $K$ be a compact subset of $X$. Recall that if $U \in
\op(X)$ then $\Gamma(U;\imin \rho F) \simeq \lpro {V
\subset\subset U} \Gamma(V;F)$, where $V \in \op(X_{sa})$. We have
the chain of isomorphisms
\begin{eqnarray*}
\Gamma(K;\imin \rho F) & \simeq & \lind U \Gamma(U;\imin \rho F) \\
& \simeq & \lind U \lpro {V \subset\subset U} \Gamma(V;F) \\
& \simeq & \lind U \Gamma(U;F)
\end{eqnarray*}
where $U$ ranges through the family of subanalytic relatively
compact open subsets of $X$ containing $K$ and $V \in
\op(X_{sa})$.

Since $F$ is quasi-injective and filtrant inductive limits are
exact, the morphism $\Gamma(X;\imin \rho F) \simeq \Gamma(X;F) \to
\lind U \Gamma(U;F) \simeq \Gamma(K;\imin \rho F)$, where $U$
ranges through the family of subanalytic open subsets of $X$
containing $K$, is surjective. \qed \\

Let us consider the following subcategory of $\mod(k_{X_{sa}})$:
\begin{eqnarray*}
\P_{X_{sa}} & := & \{G \in \mod(k_{X_{sa}});\; G \text{ is
$\Ho_{k_{X_{sa}}}(\cdot,F)$-acyclic for each $F \in {\cal
J}_{X_{sa}}$}\}.
\end{eqnarray*}

This category is generating, in fact if $\{G_j\}_j$ is a filtrant
inductive system of $\R$-constructible sheaves $\otimes_j
\rho_*G_j \in \P_{X_{sa}}$ by Corollary \ref{rcinj}. Moreover
$\P_{X_{sa}}$ is stable by $\cdot \otimes K$, where $K \in
\mod_{\rc}(k_X)$. In fact if $G \in \P_{X_{sa}}$ and $F \in {\cal
J}_{X_{sa}}$ we have
$$\Ho_{k_{X_{sa}}}(G \otimes K,F) \simeq
\Ho_{k_{X_{sa}}}(G,\ho(K,F))$$ and $\ho(K,F) \in {\cal
J}_{X_{sa}}$ by Proposition \ref{hoinj}.

\begin{teo}\label{mathcalJP} The category $\P^{op}_{X_{sa}} \times \mathcal{J}_{X_{sa}}$ is injective with respect to
the functor $\Ho_{k_{X_{sa}}}(\cdot,\cdot)$.
\end{teo}
\dim\ \ (i) Let $G \in \P_{X_{sa}}$ and consider an exact sequence
$\exs{F'}{F}{F''}$ with $F' \in {\cal J}_{X_{sa}}$. We have to
prove that the sequence
$$\exs{\Ho_{k_{X_{sa}}}(G,F')}{\Ho_{k_{X_{sa}}}(G,F)}{\Ho_{k_{X_{sa}}}(G,F'')}$$
is exact. Since the functor $\Ho_{k_{X_{sa}}}(G,\cdot)$ is acyclic
on quasi-injective sheaves we obtain the result.

(ii) Let $F \in {\cal J}_{X_{sa}}$, and let $\exs{G'}{G}{G''}$ be
an exact sequence on $\P_{X_{sa}}$. Since the objects of
$\P_{X_{sa}}$ are $\Ho_{k_{X_{sa}}}(\cdot,F)$-acyclic the sequence
$$\exs{\Ho_{k_{X_{sa}}}(G'',F)}{\Ho_{k_{X_{sa}}}(G,F)}{\Ho_{k_{X_{sa}}}(G',F)}$$
is exact. \qed


\begin{cor}  The category $\P^{op}_{X_{sa}} \times \mathcal{J}_{X_{sa}}$ is injective with respect to
the functor $\ho(\cdot,\cdot)$.
\end{cor}
\dim\ \ Let us show that $\P^{op}_{X_{sa}} \times
\mathcal{J}_{X_{sa}}$ is injective with respect to the functor
$\ho(\cdot,\cdot)$. Let $G \in \P_{X_{sa}}$, and let
$\exs{F'}{F}{F''}$ be an exact sequence with $F',F,F'' \in
\mathcal{J}_{X_{sa}}$. We shall show that for each $U \in
\op(X_{sa})$ the sequence
$$\exs{\Gamma(U;\ho(G,F''))}{\Gamma(U;\ho(G,F))}{\Gamma(U;\ho(G,F'))}$$
is exact. This is equivalent to show that for each $U \in
\op(X_{sa})$ the sequence
$$\exs{\Ho_{k_{X_{sa}}}(G_U,F'')}{\Ho_{k_{X_{sa}}}(G_U,F)}{\Ho_{k_{X_{sa}}}(G_U,F')}$$
is exact. This follows since $G_U \in \P_{X_{sa}}$. The proof of
the exactness in $\P^{op}_{X_{sa}}$ is similar. \qed


\subsection{The functor $\rho_!$}

We have seen that the functor $\imin \rho:\mod(k_{X_{sa}}) \to
\mod(k_X)$ has a left adjoint $\rho_!:\mod(k_X) \to
\mod(k_{X_{sa}})$. The functor $\rho_!$ is fully faithful and
exact. In particular, for $U \in \op(X)$ one has $\rho_!k_U \simeq
\lind {V \subset \subset U} \rho_*k_V$, where $V \in \op_{sa}(X)$.

\begin{prop}\label{roS} Let $S$ be a closed subset of $X$. Then
$\rho_!k_S \simeq \lind {W \supset S} \rho_*k_{\overline{W}}$,
where $W \in \op_{sa}(X)$.
\end{prop}
\dim\ \ (i) Let $U=X \setminus S$. Since $\rho_!$ is exact we have
an exact sequence
$$\exs{\rho_!k_U}{\rho_!k_X}{\rho_!k_S}.$$
On the other hand, let $V \in \op^c_{sa}(X)$ and $V \subset\subset
U$. We have an exact sequence $\exs{k_V}{k_X}{k_{X \setminus V}}.$
Since $\rho_*$ is exact on $\mod_{\rc}(k_X)$ the sequence
$\exs{\rho_*k_V}{\rho_*k_X}{\rho_*k_{X \setminus V}}$ is exact.
Applying the exact $\lind {V \subset\subset U}$ we obtain an exact
sequence
$$\exs{\lind {V \subset\subset U}\rho_*k_V}{\rho_*k_X}{\lind {V \subset\subset U}\rho_*k_{X \setminus V}}.$$
We have $\lind {V \subset\subset U}\rho_*k_V \simeq \rho_!k_U$ and
$\rho_*k_X \simeq \rho_!k_X$. Hence $\rho_!k_S \simeq \lind {V
\subset\subset U}\rho_*k_{X \setminus V}.$

(ii) We shall show that for each $U' \in \op^c_{sa}(X)$ the
natural morphism
\begin{equation}\label{VUWS}
\lind {V \subset\subset U}\Gamma(U';k_{X \setminus V}) \to \lind
{W \supset S} \Gamma(U';k_{\overline{W}})
\end{equation}
is an isomorphism.
We shall see that  for each $W \in \op_{sa}(X)$ with $W \supset S$
there exists $W' \in \op_{sa}(X)$ such that $X \setminus
\overline{W'} \subset\subset U$ and $\overline{W} \cap
U'=\overline{W'} \cap U'$. Set $W'=W \cup (X \setminus
\overline{U'})$. Since $U'$ is relatively compact, $X \setminus
\overline{W'} \subset\subset U$, and $\overline{W} \cap
U'=\overline{W'} \cap U'$ by construction. Then $$\lind {V
\subset\subset U}\Gamma(U';k_{X \setminus V}) \simeq \lind {(X
\setminus \overline{W}) \subset\subset U}
\Gamma(U';k_{\overline{W}}) \simeq \lind {W \supset S}
\Gamma(U';k_{\overline{W}}).\ \ $$\qed


\begin{nt} Let $Z=U \cap S$, where $U \in \op(X)$ and let $S$ be a
closed subset of $X$. Let $F \in \mod(k_{X_{sa}})$. We set for
short ${}_ZF=F \otimes \rho_!k_Z$
\end{nt}

\begin{lem} Let $F \in \mod(k_{X_{sa}})$. Let $U \in \op(X)$ and let $S$ be a closed
subset of $X$.
\begin{itemize}
\item[(i)] One has ${}_UF \simeq \lind {V \subset\subset U}F_V \simeq \lind {V \subset\subset U}
\Gamma_{\overline{V}}F$, $V \in \op^c_{sa}(X)$.
\item[(ii)] One has ${}_SF \simeq \lind {W \supset S}F_{\overline{W}} \simeq \lind {W \supset S}
\Gamma_WF$, $W \in \op_{sa}(X)$.
\end{itemize}
\end{lem}
\dim\ \ (i) The first isomorphism is obvious. Let us show the
second isomorphism. We have the chain of isomorphisms
$$\lind {V \subset\subset U}F_V \osi \lind {V,V' \subset\subset U}(\Gamma_{\overline{V'}}F)_V \iso \lind {V' \subset\subset U}
\Gamma_{\overline{V'}}F,$$ where $V,V'$ range through the family
of subanalytic open subsets of $X$.

 The proof of (ii) is
similar.
\qed \\

\begin{prop}\label{ZhoZ} Let $Z$ be a locally closed
subset of $X$. Let $G \in \mod_{\rc}(k_X)$ and $F \in
\mod(k_{X_{sa}})$. Then ${}_Z\ho(G,F) \simeq \ho(G,{}_ZF)$.
\end{prop}
\dim\ \ (i) Let $U \in \op(X)$. For each $U' \in \op^c(X_{sa})$ we
have the chain of isomorphisms
\begin{eqnarray*}
\Gamma(U';{}_U\ho(G,F)) & \simeq & \Ho_{k_{X_{sa}}}(k_{U'},\lind {V \subset\subset U}\Gamma_{\overline{V}}\ho(G,F)) \\
& \simeq & \lind {V \subset\subset U}\Ho_{k_{X_{sa}}}(k_{U'},\Gamma_{\overline{V}}\ho(G,F)) \\
& \simeq & \lind {V \subset\subset U}\Ho_{k_{X_{sa}}}(k_{U' \cap \overline{V}},\ho(G,F)) \\
& \simeq & \lind {V \subset\subset U}\Ho_{k_{X_{sa}}}(G_{U' \cap \overline{V}},F) \\
& \simeq & \lind {V \subset\subset U}\Ho_{k_{X_{sa}}}(G_{U'}, \Gamma_{\overline{V}}F) \\
& \simeq & \lind {V \subset\subset U}\Ho_{k_{X_{sa}}}(k_{U'}, \ho(G,\Gamma_{\overline{V}}F)) \\
& \simeq & \Ho_{k_{X_{sa}}}(k_{U'}, \lind {V \subset\subset U}\ho(G,\Gamma_{\overline{V}}F)) \\
& \simeq & \Ho_{k_{X_{sa}}}(k_{U'}, \ho(G,\lind {V \subset\subset U}\Gamma_{\overline{V}}F)) \\
& \simeq & \Gamma(U',\ho(G,{}_UF)),
\end{eqnarray*}
where $V \in \op_{sa}(X)$.

(ii) If $S$ is a closed subset of $X$ the proof is similar.
\qed \\


\begin{prop}\label{qinjrho!} Let $F \in \mod(k_{X_{sa}})$ be quasi-injective. Then
 $\rho_!K \otimes F$ is quasi-injective for each $K \in
 \mod(k_X)$.
\end{prop}
\dim\ \ (i) Let us show the result when $K=k_Z$, for a locally
closed subset $Z$ of $X$. Let $G \in \mod_{\rc}^c(k_X)$. We have
\begin{eqnarray*}
\Ho_{k_{X_{sa}}}(G,{}_ZF) & \simeq & \Gamma(X; \ho(G,{}_ZF) ) \\
& \simeq & \Gamma(X; {}_Z\ho(G,F) ) \\
& \simeq & \Gamma(X;\imin \rho {}_Z\ho(G,F)) \\
& \simeq & \Gamma(X;(\imin \rho \ho(G,F))_Z).
\end{eqnarray*}

 Since $F$ is
quasi-injective, $\ho(G,F)$ is quasi-injective. Then by
Proposition \ref{qsoft} the sheaf $(\imin \rho \ho(G,F))_Z$ is
c-soft and it is injective with respect to the functor
$\Gamma(X,\cdot)$. Hence
 the functor $\Gamma(X;\imin \rho \ho(\cdot,F)_Z)$ is exact on
$\mod_{\rc}^c(k_X)$.

(ii) Let $K \in \mod(k_X)$. There exists an epimorphism $\oplus_{i
\in I}k_{U_i} \twoheadrightarrow K$ with $U_i \in \op_{sa}(X)$ for
each $i$. Let $K_J$ be the image of $\oplus_{i \in J}k_{U_i}$,
with $J \subset I$ finite. We have $K \simeq \lind J K_J$, hence
$\rho_!K \simeq \lind J \rho_!K_J$ since $\rho_!$ commutes with
$\Lind$. It is enough to prove the result for $K_J$. We argue by
induction on the cardinal of $J$. Set $K=K_J$. If $|J|=1$ then
$K\simeq k_Z$ with $Z$ locally closed subset of $X$ and the result
follows from (i).

Let us show $n-1 \Rightarrow n$. There is an epimorphism
$\oplus_{i =1}^nk_{U_i} \twoheadrightarrow K$. Let $K_1$ be the
image of $k_{U_1} \to K$ and let $K_2=K/K_1$. We have a
commutative diagram
$$
\xymatrix{ 0 \ar[r] & k_{U_1} \ar[d] \ar[r] & \oplus_{i
=1}^nk_{U_i}
\ar[d] \ar[r] & \oplus_{i =2}^nk_{U_i} \ar[d] \ar[r] & 0 \\
0 \ar[r] & K_1 \ar[r] & K \ar[r] & K_2 \ar[r] & 0, }
$$
where the vertical arrows are surjective, and the rows are exact.
By the exactness of $\rho_!$ and $\otimes$ we obtain the exact
sequence $\exs{\rho_!K_1 \otimes F}{\rho_!K \otimes F}{\rho_!K_2
\otimes F}$ is exact. By the inductive hypothesis $\rho_!K_1
\otimes F$ and $\rho_!K_2 \otimes F$ are quasi-injective, then
$\rho_!K \otimes F$ is quasi-injective.
\qed \\

\begin{prop}\label{horho!} Let $F \in \mod(k_{X_{sa}})$, $G \in \mod_{\rc}(k_X)$
and let $K \in \mod(k_X)$. One has the isomorphism $\ho(G,F)
\otimes \rho_!K \simeq \ho(G,F \otimes \rho_!K)$.
\end{prop}
\dim\ \ Both sides are left exact with respect to $F$. Hence we
may assume that $F$ is quasi-injective. Since quasi-injective
sheaves are $\ho(G,\cdot)$-injective, both sides are exact with
respect to $K$. Moreover as a consequence of Proposition
\ref{hosa} both sides commute with filtrant $\Lind$ with respect
to $K$. We may reduce to the case $K=k_U$, with $U \in
\op_{sa}(X)$. Then the result follows from Proposition \ref{ZhoZ}.
\qed \\

\section{Derived category}\label{2}

As usual, we denote $D(k_{X_{sa}})$ the derived category of
$\mod(k_{X_{sa}})$ and its full subcategory consisting of bounded
(resp. bounded below, resp. bounded above) complexes is denoted by
$D^b(k_{X_{sa}})$ (resp. $D^+(k_{X_{sa}})$, resp.
$D^-(k_{X_{sa}})$).

\subsection{The category $D^b_{\rc}(k_{X_{sa}})$}

 As usual we denote by $D^b_{\rc}(k_{X_{sa}})$ (resp.
$D^b_{\rc}(k_{X_{sa}})$) the full subcategory of $D^b(k_X)$ (resp.
$D^b(k_{X_{sa}})$) consisting of objects with $\R$-constructible
cohomology.\\

Recall that $\rho:X \to X_{sa}$ is the natural morphism of sites.
It induces the functor $\rho_*:\mod(k_X) \to \mod(k_{X_{sa}})$.

\begin{lem} Let $F \in \mod_{\rc}(k_X)$. Then $R^j\rho_* F=0$
for each $j \neq 0$.
\end{lem}
\dim\ \ The sheaf $R^j\rho_*F$ is the sheaf associated to the
presheaf $V \to R^j\Gamma(V;F)$. We have to show that
$R^j\Gamma(V;F)=0$ for $j \neq 0$ on a family of generators of the
topology of $X_{sa}$. This means that for each $V \in
\op^c(X_{sa})$ and for each $j \neq 0$, there exists $I$ finite
and $\{V_i\}_{i \in I} \in \cov(V_{sa})$ such that
$R^j\Gamma(V_i;R\rho_*F) \simeq R^j\Gamma(V_i;F)=0$.

We use the notation of \cite{KS90}. There exists a locally finite
stratification $\{X_i\}_{i \in I}$ of $X$ consisting of
subanalytic subsets such that for all $j \in \Z$ and all $i \in I$
the sheaf $F|_{X_i}$ is locally constant. By the triangulation
theorem there exist a simplicial complex $(S,\Delta)$ and a
subanalytic homeomorphism $\psi:|S| \iso X$ compatible with the
stratification and such that $V$ is a finite union of the images
by $\psi$ of open subsets $V(\sigma)$ of $|S|$, where
$V(\sigma)=\bigcup_{\tau \in \Delta, \tau \supset \sigma}|\tau|$.
By Proposition 8.1.4 of \cite{KS90} we have
$R^j\Gamma(\psi(V(\sigma));F)=0$ for each $\sigma$ and for each $j
\neq 0$. The result follows because $V=\bigcup_{\psi(|\sigma|)
\subset V} \psi(V(\sigma))$. \qed \\

Since $\R$-constructible sheaves are injective with respect to the
functor $\rho_*$,
the following diagram of derived categories is quasi-commutative.
\begin{equation}\label{Drhorc}
\xymatrix{D^b_{\rc}(k_X)  \ar@ <2pt> [rr]^{R\rho_*} &&
 D^b_{\rc}(k_{X_{sa}}) \ar@ <2pt> [ll]^{\imin \rho} \\
D^b(\mod_{\rc}(k_X)) \ar[u]^\wr \ar[urr]_{\rho_*} && }
\end{equation}

\begin{teo}\label{eqrhorc} One has the equivalence of categories
$$D^b_{\rc}(k_X) \simeq D^b(\mod_{\rc}(k_X)) \simeq D^b_{\rc}(k_{X_{sa}}) .$$
\end{teo}
\dim\ \ By dévissage, to prove the equivalence between
$D^b(\mod_{\rc}(k_X))$ and $D^b_{\rc}(k_{X_{sa}})$ it is enough to
check that the functor $\rho_*$ in \eqref{Drhorc} is fully
faithful. We have $\imin \rho \circ \rho_* \simeq \id$ and the
result follows.


The equivalence between $D^b(\mod_{\rc}(k_X))$ and
$D^b_{\rc}(k_X)$ was shown by Kashiwara in \cite{Ka84}. \qed \\

\subsection{Operations in the derived category}

Let us study the operations in the derived category of
$\mod(k_{X_{sa}})$. Let $f:X \to Y$ be an analytic map. Since
$\mod(k_{X_{sa}})$ has enough injectives, then the derived
functors
\begin{eqnarray*}
\rh & : &  D^-(k_{X_{sa}})^{op} \times D^+(k_{X_{sa}}) \to D^+(k_{X_{sa}}),\\
Rf_* & : &  D^+(k_{X_{sa}}) \to D^+(k_{Y_{sa}}),\\
Rf_{!!} & : &  D^+(k_{X_{sa}}) \to D^+(k_{Y_{sa}}),
\end{eqnarray*} are well defined.

\begin{prop} Let $f:X \to Y$ be an analytic map. Then
\begin{itemize}
\item[(i)] The functors $Rf_*$ and $\rh$ commute with $R\rho_*$.
\item[(ii)] The functors $Rf_*$ and $Rf_{!!}$ commute with $\imin \rho$.
\item[(iii)] We have $R(g \circ f)_* \simeq Rg_* \circ Rf_*\ \ \text{and}\ \
R(g \circ f)_{!!} \simeq Rg_{!!} \circ Rf_{!!}.$
\item[(iv)] The functor $R^kf_{!!}:\mod(k_{X_{sa}}) \to \mod(k_{Y_{sa}})$ commutes with small filtrant inductive limits for each $k
\in \Z$.
\item[(v)] If $F \in D^+(k_{X_{sa}})$ and $f$ is proper on
$\supp (F)$, then $Rf_{!!} \simeq Rf_*$.
\end{itemize}
\end{prop}
\dim\ \ (i) The functor $\rho_*$ sends injective sheaves to
injective sheaves, then $Rf_*$ and $\rh$ commute with $R\rho_*$.

(ii) Since $\imin \rho$ has a left adjoint it sends injective
sheaves to injective sheaves. Then $Rf_*$ and $Rf_{!!}$ commute
with $\imin \rho$.

(iii) The functor $f_*$ (resp. $f_{!!}$) sends injective sheaves
to injective (resp. quasi-injective) sheaves. Then $R(g \circ f)_*
\simeq Rg_* \circ Rf_*\ \ \text{and}\ \ R(g \circ f)_{!!} \simeq
Rg_{!!} \circ Rf_{!!}.$

(iv) Quasi-injective objects of $\mod(k_{X_{sa}})$ are stable by
filtrant $\Lind$, and the functor $f_{!!}$ commutes with such
limits. Then $R^kf_{!!}$ commutes with filtrant $\Lind$ for each
$k \in \Z$.

(v) We can find a representative $F'$ of $F$ in
$K^+(\mathcal{J}_{X_{sa}})$ with $f$ proper on $\supp(F')$. Then
the result follows from the non derived case. \qed \\

\begin{prop}\label{Rhosa} Let $F = \lind i F_i$ with $F_i \in
\mod(k_{X_{sa}})$ and let $G \in D^b_{\rc}(k_X)$.
One has $R^k\ho(G,F) \simeq \lind i R^k\ho(G,F_i)$ for each $k \in
\Z$.
\end{prop}
\dim\ \  There exists (see \cite{KS}, Corollary 9.6.7) an
inductive system of injective resolutions $I^\bullet_i$ of $F_i$.
Then
$\lind i I_i^\bullet$ is a complex of quasi-injective objects
quasi-isomorphic to $F$. Each object of
$(\mod_{\rc}(k_X)^{op},\mathcal{J}_{X_{sa}})$ is
$\ho(\cdot,\cdot)$-acyclic.
Proposition \ref{hosa} implies  the isomorphism
$$\ho(G,\lind i I_i^\bullet) \simeq  \lind i
\ho(G,I^\bullet_i)$$
and the result follows.
\qed \\


\begin{prop}\label{Rhorho!} Let $F \in D^+(k_{X_{sa}})$, $G \in D^b_{\rc}(k_X)$
and let $K \in D^+(k_X)$. One has the isomorphism $\rh(G,F)
\otimes \rho_!K \simeq \rh(G,F \otimes \rho_!K)$.
\end{prop}
\dim\ \ Let $I^\bullet$ be a quasi-injective resolution of $F$. By
Proposition \ref{qinjrho!} we have that $I^\bullet \otimes
\rho_!K$ is a complex of quasi-injective objects. Each object of
$(\mod_{\rc}(k_X)^{op},\mathcal{J}_{X_{sa}})$ is
$\ho(\cdot,\cdot)$-acyclic. Hence we are reduced to prove the
isomorphism $\ho(G,I^\bullet) \otimes \rho_!K \simeq
\ho(G,I^\bullet \otimes \rho_!K)$. The result follows from
Proposition \ref{horho!}. \qed \\

\begin{prop} Let $U \in \op^c(X_{sa})$. Let $F \in \mod(k_{X_{sa}})$ be
quasi-injective. Then $F_U$ is $\Gamma(V,\cdot)$-acyclic for each
$V \in \op(X_{sa})$.
\end{prop}
\dim\ \ Since $F_U$ has compact support, we may suppose that $V$
is relatively compact. Let $S=X \setminus U$. Since $F$ is
quasi-injective and filtrant $\Lind$ are exact, the morphism
$\Gamma(V;F) \to \lind {W \supset S \cap V}\Gamma(W; F) \iso
\Gamma(V;F_S)$ is surjective. Consider the exact sequence
$\exs{F_U}{F}{F_S}$. We get the exact sequence
$$\exs{\Gamma(V;F_U)}{\Gamma(V;F)}{\Gamma(V;F_S)}.$$
By Proposition \ref{FSqinj} $F$ and $F_S$ are quasi-injective,
hence $\Gamma(V;\cdot)$-acyclic. This implies that $F_U$ is
$\Gamma(V;\cdot)$-acyclic. \qed \\

\begin{cor}\label{Uf!!ac} Let $f:X \to Y$ be a real analytic map and let
$U \in \op^c(X_{sa})$. Let $F \in \mod(k_{X_{sa}})$ be
quasi-injective. Then $F_U$ is $f_{!!}$-acyclic.
\end{cor}
\dim\ \ Since $F_U$ has compact support, $Rf_{!!}F_U \simeq
Rf_*F_U$. The result follows because $F_U$ is $\Gamma(\imin
f(V);\cdot)$-acyclic for each $V \in \op(Y_{sa})$. \qed \\

\begin{lem}\label{rcf!!ac} Let $F$ be
quasi-injective object of $\mod(k_{X_{sa}})$ and let $G \in
\mod^c_{\rc}(k_X)$. Then $F \otimes \rho_*G$ is $f_{!!}$-acyclic.
\end{lem}
\dim\ \ Let $G \in \mod_{\rc}^c(k_X)$. Then $G$ has a resolution
$$0 \to \oplus_{i_1 \in I_1} k_{U_{i_1}} \to \ldots \to \oplus_{i_n \in I_n} k_{U_{i_n}}
\stackrel{\varphi}{\to} G \to 0.$$ Where $I_j$ is finite and
$U_{i_j} \in \op^c_{sa}(X)$ for each $i_j \in I_j$, $j \in
\{1,\ldots ,n\}$. Let us argue by induction on the length $n$ of
the resolution.

If $n=1$, then $G$ is isomorphic to a finite sum $\oplus_i
k_{U_i}$, with $U_i \in \op^c_{sa}(X)$, and the result follows
from Corollary \ref{Uf!!ac}.

Let us show $n-1 \Rightarrow n$. The sequence $\exs{\ker
\varphi}{\oplus_{i_n \in I_n} k_{U_{i_n}}}{G}$ is exact. The sheaf
$\ker \varphi$ belongs to $\mod_{\rc}^c(k_X)$ and it has a
resolution of length $n-1$. Applying $F \otimes \rho_*(\cdot)$ we
get the exact sequence
$$\exs{F \otimes \rho_*\ker \varphi}{\oplus_{i_n \in I_n} F_{U_{i_n}}}{F \otimes \rho_*G}.$$
By the induction hypothesis $F \otimes \rho_*\ker \varphi$ is
$f_{!!}$-acyclic. Moreover $\oplus_{i_n \in I_n}F_{U_{i_n}}$ is
$f_{!!}$-acyclic, then $F \otimes \rho_*G$ is $f_{!!}$-acyclic.
\qed \\

\begin{prop}\label{saf!!ac} Let $F$ be
quasi-injective object of $\mod(k_{X_{sa}})$ and let $G \in
\mod(k_{X_{sa}})$. Then $F \otimes G$ is $f_{!!}$-acyclic.
\end{prop}
\dim\ \ Let $G \simeq \lind i \rho_*G_i$ with $G_i \in
\mod^c_{\rc}(k_X)$ for each $i$. Since the functors $\otimes$ and
$R^kf_{!!}$ commute with filtrant $\Lind$ we have
$$R^kf_{!!}(F \otimes \lind i \rho_*G_i) \simeq \lind i R^kf_{!!}(F
\otimes \rho_*G_i)=0$$ if $k \neq 0$ by Lemma \ref{rcf!!ac}.
\qed \\

\begin{prop}\label{Rbcsa} Let $F \in D^+(k_{X_{sa}})$ and $G \in
D^+(k_{Y_{sa}})$. Then
$$Rf_{!!}F \otimes G \simeq Rf_{!!}(F \otimes \imin f G).$$
\end{prop}
\dim\ First assume that $F \in \mod(k_{X_{sa}})$ is injective. By
Proposition \ref{saf!!ac} $F \otimes \imin fG$ is
$f_{!!}$-acyclic.

Now let $F \in D^+(k_{X_{sa}})$ and $G \in D^+(k_{Y_{sa}})$. Let
$F'$ be a complex of injective sheaves quasi-isomorphic to $F$.
Then
$$Rf_{!!}F \otimes G \simeq f_{!!}F' \otimes G \simeq f_{!!}(F' \otimes \imin f G)
\simeq Rf_{!!}(F \otimes \imin f G),$$ where the second
isomorphism follows from Proposition \ref{bcsa}. \qed \\

Now let us consider a cartesian square
$$
\xymatrix{X'_{sa} \ar[r]^{f'} \ar[d]^{g'} & Y'_{sa} \ar[d]^g \\
X_{sa} \ar[r]^f &Y_{sa}}
$$
\begin{prop} Let $F \in D^+(k_{X_{sa}})$. Then $\imin g Rf_{!!}F
\simeq Rf'_{!!}g'{}^{-1}F$.
\end{prop}
\dim\ \ We have an isomorphism $f'_{!!}g'{}^{-1} \simeq \imin g
f_{!!}$, and $R(\imin g f_{!!}) \simeq \imin g Rf_{!!}$ since
$\imin g$ is exact. Then we obtain a morphism $\imin g Rf_{!!} \to
Rf'_{!!}g'{}^{-1}$. It is enough to prove that for any $k \in \Z$
and for any $F \in \mod(k_{X_{sa}})$ we have $\imin g R^kf_{!!}F
\iso R^kf'_{!!}g'{}^{-1}F$. Since both sides commute with filtrant
$\Lind$, we may assume $F=\rho_*G$ with $G \in \mod^c_{\rc}(k_X)$.
Moreover since $\supp (G)$ is compact, $f'$ is proper on $\supp
(g'{}^{-1}G)$. Then both sides commute with $\rho_*$ and the
result follows from the corresponding one for classical sheaves.
\qed \\

As in classical sheaf theory, the K\"unneth formula follows from
the projection formula and the base change formula.

\begin{prop} Consider a cartesian square
$$
\xymatrix{X'_{sa} \ar[r]^{f'} \ar[d]^{g'} \ar@{-->} [dr]^\delta  & Y'_{sa} \ar[d]^g \\
X_{sa} \ar[r]^f &Y_{sa}}
$$
where $\delta=fg'=gf'$. There is a natural isomorphism
$$R\delta_{!!}(g'{}^{-1}F \otimes f'^{-1}G) \simeq Rf_{!!}F \otimes Rg_{!!}G$$
for $F \in D^+(k_{X_{sa}})$ and $G \in D^+(k_{Y'_{sa}})$.
\end{prop}
\dim\ \ Using the projection formula and the base change formula
we deduce
$$Rf'_{!!}(g'{}^{-1}F \otimes f'^{-1}G) \simeq Rf'_{!!}g'^{-1}F \otimes
G \simeq \imin g Rf_{!!}F \otimes G.$$

Using the projection formula once again we find
$$R\delta_{!!}(g'{}^{-1}F \otimes f'^{-1}G) \simeq Rg_{!!}Rf'_{!!}(g'{}^{-1}F \otimes
f'^{-1}G) \simeq Rf_{!!}F \otimes Rg_{!!}G$$
and the result follows. \qed \\


\begin{prop}\label{nosheaf} Let $G \in D^b_{\rc}(k_Y)$ and let $F
\in D^+(k_{X_{sa}})$. Then the natural morphism
$$Rf_{!!}\rh(\imin f G,F) \to \rh(G,Rf_{!!}F)$$
is an isomorphism.
\end{prop}
\dim\ \ The morphism is obtained as in the non derived case. Let
us show that it is an isomorphism. Let $F'$ be a complex of
injective sheaves quasi-isomorphic to $F$. Then
\begin{eqnarray*}
Rf_{!!}\ho(\imin f G,F) & \simeq & f_{!!}\ho(\imin f G,F') \\
& \simeq & \ho(G,f_{!!}F') \\
& \simeq & \rh(G,Rf_{!!}F),
\end{eqnarray*}
 where the second
isomorphism follows from Proposition \ref{nomod}. \qed \\

\subsection{Vanishing theorems on $\mod(k_{X_{sa}})$}

In this Section we give some results on the vanishing of the
cohomology of sheaves on a subanalytic site.

\begin{df} The quasi-injective dimension of the category $\mod(k_{X_{sa}})$ is
the smallest $n \in \N \cup \{\infty\}$ such that for any $F \in
\mod(k_{X_{sa}})$ there exists an exact sequence
$$0 \to F \to \I^0 \to \cdots \to I^n \to 0$$
with $I^j$ quasi-injective for $0 \leq j \leq n$.
\end{df}

\begin{prop}\label{fqid} The category $\mod(k_{X_{sa}})$ has finite
quasi-injective dimension.
\end{prop}
\dim\ \ Let $F \in \mod(k_{X_{sa}})$. Then $F=\lind i \rho_*F_i$,
with $F_i \in \mod^c_{\rc}(k_X)$.  There exists (see \cite{KS},
Corollary 9.6.7) an inductive system of injective resolutions
$I^\bullet_i$ of $F_i$. By Proposition 3.3.11 of \cite{KS90}, the
category $\mod(k_X)$ has finite homological dimension. Then we may
assume that $I_i^\bullet$ has length $N_0<\infty$ for each $i$.
Since $F_i$ is $\rho_*$-injective for each $i$,
$\rho_*I_i^\bullet$ is an injective resolution of $\rho_*F_i$ of
length $N_0$. Taking the inductive limit we find that $\lind i
\rho_*I_i^\bullet$ is a resolution of $F$ of length $N_0$, and
$\lind i\rho_*I_i^j \in \mathcal{J}_{X_{sa}}$ for each $j$.
\qed \\

\begin{cor} Let $f:X \to Y$ be a real analytic map, and let $F \in \mod_{\rc}(k_X)$. The functors $f_*$,
$f_{!!}$ and $\ho(F,\cdot)$ have finite cohomological dimension.
\end{cor}

\begin{prop} Let $F \in \mod(k_X)$ and let $G \in
\mod(k_{X_{sa}})$. There exists a finite $j_0 \in \N$ such that
$$R^j\ho(\rho_!F,G)=0 \ \ \text{for $j>j_0.$}$$
\end{prop}
\dim\ \ Let $U \in \op(X_{sa})$. We have the chain of isomorphisms
\begin{eqnarray*}
R\Gamma(U;\rh(\rho_!F,G)) & \simeq & \Rh_{k_{X_{sa}}}(\rho_!F,R\Gamma_UG) \\
& \simeq & \Rh_{k_X}(F,\imin \rho R\Gamma_UG).
\end{eqnarray*}
The functor $R\Gamma_U$ has finite cohomological dimension, and
the homological dimension of the category $\mod(k_X)$ is finite.
Hence we can find a finite $j_0 \in \N$ such that
$R^j\Gamma(U;\rh(\rho_!F,G))$ vanishes for $j>j_0$ and for each $U
\in \op(X_{sa})$. This shows the result. \qed \\

\begin{oss} We have seen that the functor $\ho(F,\cdot)$ has
finite cohomological dimension when $F$ is $\R$-constructible and
when $F=\rho_!G$ with $G \in \mod(k_X)$. We do not know if the
cohomological dimension is finite for any $F \in
\mod(k_{X_{sa}})$. Indeed we do not know if the homological
dimension of $\mod(k_{X_{sa}})$ is finite or not.
\end{oss}

\subsection{Duality}

In the following we find a right adjoint to the functor $Rf_{!!}$,
denoted by $f^!$, and we calculate it by decomposing $f$ as the
composite of a closed embedding and a submersion.\\

The subcategory $\mathcal{J}_{X_{sa}}$ of quasi-injective objects
and the functor $f_{!!}$ have the following properties:
\begin{equation}\label{Brownhyp}
  \begin{cases}
    \text{(i) $\mathcal{J}_{X_{sa}}$ is cogenerating}, \\
    \text{(ii) $\mod(k_{X_{sa}})$ has finite quasi-injective dimension}, \\
    \text{(iii) $\mathcal{J}_{X_{sa}}$ is $f_{!!}$-injective}, \\
    \text{(iv) $\mathcal{J}_{X_{sa}}$ is closed by small $\oplus$}, \\
    \text{(v) $f_{!!}$ commutes with small $\oplus$}.
  \end{cases}
  \end{equation}
As a consequence of the Brown representability theorem (see
\cite{KS}, Corollary 14.3.7 for details) we find a right adjoint
to the functor $Rf_{!!}$.

\begin{teo}\label{f!D+} (i) The functor $Rf_{!!}:D(k_{X_{sa}}) \to D(k_{Y_{sa}})$
admits a right adjoint. We denote by $f^!:D(k_{Y_{sa}}) \to
D(k_{X_{sa}})$ the adjoint functor.

(ii) Let $G \in D^+(k_{Y_{sa}})$. Then $f^!G \in D^+(k_{X_{sa}})$.
\end{teo}

\begin{oss} As in classical sheaf theory, one can prove by
adjunction the dual projection formula and the dual base change
formula.
\end{oss}


\begin{prop}\label{f!lim} The functor $f^!$ commutes
with $R\rho_*$, and the functor $H^kf^!:\mod(k_{Y_{sa}}) \to
\mod(k_{X_{sa}})$ commutes with filtrant $\Lind$.
\end{prop}
\dim\ \ Since $Rf_{!!}$ commutes with $\imin \rho$, then $f^!$
commutes with $R\rho_*$ by adjunction.

Let us show that $H^kf^!$ commutes with $\Lind$. Let $\{F_i\}_i$
be a filtrant inductive system in $\mod(k_{Y_{sa}})$. Remark that
$\lind i H^kf^!F_i$ (resp. $H^kf^!\lind iF_i$) is the sheaf
associated to the presheaf $U \to \lind i R^k\Gamma(U;f^!F_i)$
(resp. $U \to  R^k\Gamma(U;f^!\lind iF_i)$), for $U \in
\op^c(X_{sa})$.

We will show the isomorphism $R^k\Gamma(U;f^!\lind i F_i)
\stackrel{\sim}{\gets} \lind i R^k\Gamma(U;f^!F_i)$ for each $U
\in \op^c_{sa}(X)$. By adjunction it is enough to prove the
isomorphism $R^k\Ho_{k_{Y_{sa}}}(Rf_!k_U,\lind i F_i) \simeq \lind
i R^k\Ho_{k_{Y_{sa}}}(Rf_!k_U,F_i).$

Let $\mathcal{J}_{Y_{sa}}$ be the fa\-mi\-ly of quasi-injective
objects of $\mod(k_{Y_{sa}})$. Each object of
$(\mod^c_{\rc}(k_Y)^{op},\mathcal{J}_{Y_{sa}})$ is
$\Ho_{k_{Y_{sa}}}(\cdot,\cdot)$-acyclic. Moreover
$\mathcal{J}_{Y_{sa}}$ is stable by filtrant inductive limits.
There exists (see \cite{KS}, Corollary 9.6.7) an inductive system
of injective resolutions $I^\bullet_i$ of $F_i$. Then $\lind
iI_i^\bullet$ is a quasi-injective resolution of $\lind i F_i$. We
have
$$\Ho_{K^+(k_{Y_{sa}})}(Rf_!k_U,\lind i I^\bullet_i) \simeq \lind i \Ho_{K^+(k_{Y_{sa}})}(Rf_!k_U,I^\bullet_i)$$
and the result follows. \qed \\

\begin{cor} Let $F \in D^b(k_{Y_{sa}})$. Then $f^!F \in D^b(k_{X_{sa}})$.
\end{cor}
\dim\ \ We may reduce to the case $F \in \mod(k_{Y_{sa}})$. Then
$F \simeq \lind i \rho_*F_i$ with $F_i \in \mod_{\rc}(k_Y)$ for
each $i$. By Proposition \ref{f!lim} we have
$$H^kf^!F \simeq H^kf^!\lind i\rho_*F_i \simeq \lind i\rho_*
H^kf^!F_i,$$ and $H^kf^!F_i=0$ if $k>j_0$ for a fixed $j_0\in \N$
and for each $i$. \qed \\

\begin{prop} Let $F \in D^+(k_{Y_{sa}})$ and let $G \in D^+(k_Y)$. Then
one has the isomorphism $f^!(F \otimes \rho_!G) \simeq f^!F
\otimes \rho_!\imin f G.$
\end{prop}
\dim\ \ We have the chain of morphisms
$$Rf_{!!}(f^!F \otimes \rho_!\imin f G) \simeq Rf_{!!}f^!F \otimes \rho_!G \to F \otimes \rho_!G,$$
by adjunction we obtain the desired morphism. To prove that it is
an isomorphism it is enough to show $R^k\Gamma(U;f^!(F \otimes
\rho_!G)) \simeq R^k\Gamma(U;f^!F \otimes \rho_!\imin f G)$ for
each $U \in \op^c(X_{sa})$ and each $k \in \Z$. We have the chain
of isomorphisms
\begin{eqnarray*}
R^k\Gamma(U;f^!(F \otimes \rho_!G)) & \simeq &
R^k\Ho_{k_{Y_{sa}}}(Rf_{!!}k_U,F \otimes \rho_!G) \\
& \simeq &
R^k\Ho_{k_{Y_{sa}}}(k_Y,\rh(Rf_{!!}k_U,F \otimes \rho_!G)) \\
& \simeq &
R^k\Ho_{k_{Y_{sa}}}(k_Y,\rh(Rf_{!!}k_U,F) \otimes \rho_!G) \\
& \simeq &
R^k\Ho_{k_{Y_{sa}}}(k_Y,Rf_{!!}\rh(k_U,f^!F) \otimes \rho_!G) \\
& \simeq &
R^k\Ho_{k_{Y_{sa}}}(k_Y,Rf_{!!}(\rh(k_U,f^!F) \otimes \imin f\rho_!G)) \\
& \simeq &
R^k\Ho_{k_{Y_{sa}}}(k_Y,Rf_{!!}(\rh(k_U,f^!F \otimes \rho_!\imin f G)) \\
& \simeq & R^k\Ho_{k_{X_{sa}}}(k_U,f^!F \otimes \rho_!\imin f G)).
\end{eqnarray*}
Here the fourth and the last isomorphism follow from the fact that
since $k_U$ has compact support, then $\rh(k_U,K)$ has compact
support for any $K \in D^+(k_{X_{sa}})$ and $Rf_{!!}\rh(k_U,K)
\simeq Rf_*\rh(k_U,K)$. \qed \\

\begin{prop}\label{closemb} Let $F \in D^+(k_{Y_{sa}})$, and let $f:X \to Y$ be
a closed embedding.
Then $f^!F \simeq \imin f\rh(k_X,F)$ and $\id \iso f^!Rf_{!!}$.
\end{prop}
\dim\ \
Since $f$ is proper, then $Rf_* \simeq Rf_{!!}$. We have the
isomorphisms $Rf_*f^!F \simeq Rf_*\rh(k_X,f^!F) \simeq
 \rh(k_X,F)$. Since $\imin f Rf_* f^!F \simeq f^!F$, then $f^!F \simeq
 \imin f \rh(k_X,F)$.

 Let $F' \in D^+(k_{X_{sa}})$. We have the
 isomorphisms
 $$
 f^!Rf_*F' \simeq \imin f \rh(k_X,Rf_*F') \simeq \imin f Rf_*\rh(k_X,F') \simeq \imin f
 Rf_*F',
 $$
and $\imin fRf_*F' \simeq F'$ since $f$ is a closed embedding.
\qed \\

 Recall
that $f$ is a topological submersion (of fiber dimension $n$) if
locally on X, $f$ is isomorphic to the projection $Y \times \R^n
\to Y$.

\begin{prop}\label{topsubsa1} Assume that $f$ is a topological submersion. Then
for $F \in D^+(k_{Y_{sa}})$ one has the isomorphism $\imin f F
\otimes f^!k_Y \iso f^!F.$
\end{prop}
\dim\ \ We have the chain of morphisms
$$Rf_{!!}(\imin f F \otimes f^!k_Y) \simeq F \otimes
Rf_{!!}f^!k_Y \to F \otimes k_Y \simeq F,$$ by adjunction we
obtain the desired morphism.

Let us show that it is an isomorphism. We may reduce to the case
$F \simeq \lind i \rho_*F_i \in \mod(k_{Y_{sa}})$. We have the
chain of isomorphisms
\begin{eqnarray*}
H^k(\imin f \lind i \rho_* F_i \otimes f^!\rho_*k_Y) & \simeq &
\lind i
\rho_*H^k(\imin f F_i \otimes f^!k_Y) \\
& \simeq & \lind i \rho_*H^kf^!F_i \\
& \simeq & H^kf^!\lind i \rho_*F_i.
\end{eqnarray*}
\qed

Using these results we can calculate explicitly the functor $f^!$.
Let $f:X \to Y$ be an analytic map. We decompose it as the
composite of a closed embedding and a submersion. In fact
$$f:X \overset{j}{\hookrightarrow} X \times Y \overset{p}{\to} Y$$
where $p$ is the projection and $j$ is the graph embedding
$j(x)=(x,f(x))$. Let $F \in D^+(k_{Y_{sa}})$. Applying
Propositions \ref{closemb} and \ref{topsubsa1} we get
$$f^!F \simeq \imin j \rh(k_{j(X)},\imin p F \otimes p^!k_Y).$$

\begin{cor}\label{topsubsacor} Assume that $f$ is a topological submersion. Then:
\begin{itemize}
\item[(i)] the functor $f^!$ commutes with $\imin \rho$,
\item[(ii)] the functor $Rf_{!!}$ commutes with $\rho_!$.
\end{itemize}
\end{cor}
\dim\ \ (i) One has the chain of isomorphisms
\begin{eqnarray*}
\imin \rho (\imin f F \otimes f^!\rho_*k_Y) & \simeq & \imin \rho
\imin f F \otimes \imin \rho f^!\rho_*k_Y \\ & \simeq & \imin f
\imin \rho F
\otimes \imin \rho f^!\rho_*k_Y \\
& \simeq & \imin f \imin \rho F \otimes \imin \rho \rho_* f^!k_Y
\\
& \simeq & \imin f \imin \rho F \otimes f^!k_Y.
\end{eqnarray*}
The result follows from Proposition \ref{topsubsa1}.

(ii) The result follows by adjunction. \qed \\

\begin{prop}\label{topsubsa2} Assume that f is a topological submersion and moreover that $Rf_!f^!k_Y
\simeq k_Y$. Then for $F \in D^+(k_{Y_{sa}})$ the morphism
$Rf_*\imin fF \to F$ is an isomorphism.
\end{prop}
\dim\ \ First let us show that $Rf_{!!}f^!k_Y \simeq k_Y$. We have
the chain of isomorphisms
\begin{eqnarray*}
Rf_{!!}f^!\rho_* k_Y & \simeq & Rf_{!!}\rho_* f^!k_Y \\
& \simeq & Rf_{!!}\rho_! f^!k_Y \\
& \simeq & \rho_! Rf_! f^!k_Y \\
& \simeq & \rho_! k_Y \\
& \simeq & \rho_* k_Y,
\end{eqnarray*}
where the second isomorphism follows because $f^!k_Y$ is locally
constant and the third from Corollary \ref{topsubsacor} (ii). It
follows from Proposition \ref{topsubsa1} that $\imin f F \simeq
\rh(f^!k_Y,f^!F)$. Then we have the chain of isomorphisms
\begin{eqnarray*}
Rf_*\imin f F & \simeq & Rf_*\rh(f^!k_Y,f^!F)\\
& \simeq & \rh(Rf_{!!}f^!k_Y,F)\\
& \simeq & F.
\end{eqnarray*}
\qed

\section{Examples of applications}\label{3}

In this Section we give some example of subanalytic sheaves. Let
$X$ be a real analytic manifold, and let $X_{sa}$ be the
associated subanalytic site. We first introduce sheaves of
$\RR$-modules, where $\RR$ is a sheaf of $k$-algebras on $X_{sa}$.
Let $\D_X$ be the sheaf of finite order differential operators on
$X$. We define the $\rho_!\D_X$-modules $\ot_X$ and $\OWX$ of
tempered and Whitney holomorphic functions respectively.
References are made to \cite{KS} for an exposition on sheaves of
rings on a Grothendieck topology.

\subsection{Modules over a $k_{X_{sa}}$-algebra}

A sheaf of $k_{X_{sa}}$-algebras (or a $k_{X_{sa}}$-algebra, for
short) is an object $\RR \in \mod(k_{X_{sa}})$ such that
$\Gamma(U;\RR)$ is a $k$-algebra for each $U \in \op(X_{sa})$. The
opposite $k_{X_{sa}}$-algebra $\RR^{op}$ is defined by setting
$\Gamma(U;\RR^{op})=\Gamma(U;\RR)^{op}$ for each $U \in
\op(X_{sa})$. A sheaf of (left) $\RR$-modules is a sheaf $F$ such
that $\Gamma(U;F)$ has a structure of (left)
$\Gamma(U;\RR)$-module for each $U \in \op(X_{sa})$.
\\

Let $\RR$ be a $k_{X_{sa}}$-algebra and denote by $\mod(\RR)$ the
category of sheaves of (left) $\RR$-modules. The category
$\mod(\RR)$ is a Grothendieck category and the family
$\{\RR_U\}_{U \in \op^c(X_{sa})}$ is a small system of generators.
Moreover the forgetful functor $for:\mod(\RR) \to
\mod(k_{X_{sa}})$ is exact.\\

In this Section we shall extend some results on
$k_{X_{sa}}$-modules, by replacing $k_{X_{sa}}$ with $\RR$. Since
the formalism is similar to that we developed previously we shall
not give proofs. The functors
\begin{eqnarray*}
\ho_\RR & : &  \mod(\RR)^{op} \times \mod(\RR) \to \mod(k_{X_{sa}}),\\
\otimes_\RR & : &  \mod(\RR^{op}) \times \mod(\RR) \to
\mod(k_{X_{sa}})
\end{eqnarray*}
 are well defined. Let us summarize their
properties:
\begin{itemize}
\item the functor $\ho_\RR$ is left exact,
\item the functor $\otimes_\RR$ is right exact and commutes with $\Lind$.
\end{itemize}


Let $X,Y$ be two real analytic manifolds, and let $f:X \to Y$ be a
morphism of real analytic manifolds. Let $\RR$ be a
$k_{Y_{sa}}$-algebra. The functors $\imin f$, $f_*$ and $f_{!!}$
induce functors
\begin{eqnarray*}
\imin f & : & \mod(\RR) \to \mod(\imin f
\RR ),\\
f_* & : &  \mod(\imin f
\RR ) \to \mod(\RR),\\
f_{!!} & : &  \mod(\imin f \RR ) \to \mod(\RR).
\end{eqnarray*}
Let us summarize their properties:
\begin{itemize}
\item the functor $\imin f$ is exact and commutes with $\Lind$ and $\otimes_\RR$,
\item the functor $f_*$ is left exact and commutes with $\Lpro$,
\item $(\imin f,f_*)$ is a pair of adjoint functors,
\item the functor $f_{!!}$ is left exact and commutes with filtrant $\Lind$.
\end{itemize}

Now we consider the derived category of sheaves of $\RR$-modules.

\begin{df} An object $F \in \mod(\RR)$ is flat if the functor $\mod(\RR^{op}) \ni G
\to G \otimes_\RR F$ is exact.
\end{df}

Small direct sums and filtrant inductive limits of flat
$\RR$-modules are flat. Since the generators of $\mod(\RR)$ are
flat, then the subcategory of $\mod(\RR)$ consisting of flat
modules is generating. Thanks to flat objects we can find a left
derived functor $\otimes^L_\RR$ of the tensor product
$\otimes_\RR$.

\begin{df} An object $F \in \mod(\RR)$ is quasi-injective if the
functor $\Ho_{k_{X_{sa}}}(\cdot,F)$ is exact in
$\mod_{\rc}^c(k_X)$ or, equivalently (see Theorem 8.7.2 of
\cite{KS}) if for each $U,V \in \op^c(X_{sa})$ with $V \subset U$
the restriction morphism $\Gamma(U;F) \to \Gamma(V;F)$ is
surjective.
\end{df}

Let $X,Y$ be two real analytic manifolds, and let $f:X \to Y$ be a
real analytic map. Let $\RR$ be a $k_{Y_{sa}}$-algebra. As in
$\S$\,\ref{quinj} one can prove that quasi-injective objects are
injective with respect to the functors $f_*$ and $f_{!!}$. The
functors $Rf_*$ and $Rf_{!!}$ are well defined and projection
formula, base change formula and K\"unneth formula remain valid
for $\RR$-modules. Moreover hypothesis \eqref{Brownhyp} are
satisfied and we have

\begin{teo} The functor $Rf_{!!}:D^+(\imin f \RR) \to D^+(\RR)$
admits a right adjoint. We denote by $f^!:D^+(\RR) \to D^+(\imin f
\RR)$ the adjoint functor.
\end{teo}

\subsection{Sheaves of $\rho_!\RR$-modules}

We will consider the case where the ring is $\rho_!\RR$, where
$\RR$ is a sheaf of $k_X$-algebras. We will also assume the
following hypothesis:
$$\text{$\RR$ has finite flat dimension.}$$

The functor $\rho_!$ induces an exact functor $\mod(\RR) \to
\mod(\rho_!\RR)$ which is left adjoint to $\imin
\rho:\mod(\rho_!\RR) \to \mod(\RR)$. We will still denote by
$\rho_!$ that functor. The functor $\rho_*:\mod(\RR) \to
\mod(\rho_!\RR)$ is well defined too, in fact the morphism $\xi_F
\in \Ho_{k_X}(\RR,\mathcal{E}nd(F))$ defines a morphism in
$\Ho_{k_{X_{sa}}}(\rho_!\RR,\mathcal{E}nd(\rho_*F))$. That follows
from the chain of isomorphism
\begin{eqnarray*}
\Ho_{k_{X_{sa}}}(\rho_!\RR,\mathcal{E}nd(\rho_*F)) & \simeq &
\Ho_{k_{X_{sa}}}(\rho_!\RR,\rho_*\mathcal{E}nd(F)) \\
& \simeq & \Ho_{k_X}(\imin \rho \rho_!\RR,\mathcal{E}nd(F)) \\
& \simeq & \Ho_{k_X}(\RR,\mathcal{E}nd(F)).
\end{eqnarray*}

We briefly summarize the properties of these functors:

\begin{itemize}
\item $\imin \rho$ commutes with $\otimes^L_\RR$, $\imin f$ and
$Rf_{!!}$,
\item $R\rho_*$ commutes with $\rh_\RR$ and $Rf_*$,
\item $\rho_!$ commutes with $\otimes^L_\RR$ and $\imin f$.
\end{itemize}

Finally we recall the following result (which has been proved in
\cite{KS01})

\begin{prop}\label{6.7.1} Denote by $\widetilde{\RR}$ the presheaf
$U \to \Gamma(\overline{U};\RR)$, where $U \in \op^c(X_{sa})$.
Suppose that $F$ is a presheaf of $\widetilde{\RR}$-modules and
denote by $F^{++}$ the sheaf associated to $F$. Then $F^{++} \in
\mod(\rho_!\RR)$.
\end{prop}
\dim\ \ Let $U \in \op(X_{sa})$, and let $r \in
\Gamma(\overline{U};\RR)$. Then $r$ defines a morphism
$\Gamma(\overline{V};\RR) \otimes \Gamma(V;F) \to \Gamma(V;F)$ for
each subanalytic $V \subset U$, hence un endomorphism of
$(F^{++})|_{U_{X_{sa}}}\simeq (F|_{U_{X_{sa}}})^{++}$. This
morphism defines a morphism of sheaves $\widetilde{\RR} \to {\cal
E}nd(F^{++})$ and $\widetilde{\RR}^{++}\simeq\rho_!\RR$ by
Proposition \ref{eta!}. Then $F^{++} \in \mod(\rho_!\RR)$. \qed

\subsection{Some examples of subanalytic sheaves}

From now on, the base field is $\CC$. Let $M$ be a real analytic
manifold. One denotes by $\C^\infty_M$ and $\db_M$ the sheaves of
$\C^\infty$ functions and Schwartz's distributions respectively,
and by $\D_M$ the sheaf of finite order
differential operators with analytic coefficients. As usual, given a sheaf $F$ on $M$, we set $D'F=\rh(F,\CC_M)$.\\

In \cite{Ka84} the author defined the functor
$$\th(\cdot,\db_M):\mod_{\rc}(\CC_M)^{op} \to \mod(\D_M)$$
in the following way: let $U$ be a subanalytic open subset of $M$
and $Z=M \setminus U$. Then the sheaf $\th(\CC_U,\db_M)$ is
defined by the exact sequence
$$\exs{\Gamma_Z\db_M}{\db_M}{\th(\CC_U,\db_M)}.$$
This functor is exact and extends as a functor in the derived
category, from $D^b_{\rc}(\CC_M)$ to $D^b(\D_M)$. Moreover the
sheaf $\th(F,\db_M)$ is soft for any $\R$-constructible sheaf
$F$.\\

\begin{df} One denotes by $\dbt_M$ the presheaf of tempered
distributions on $M_{sa}$ defined as follows:
$$U \mapsto \Gamma(M;\db_M)/\Gamma_{M\setminus U}(M;\db_M).$$
\end{df}
 As a consequence of the \L ojasievicz's inequalities
\cite{Lo59}, for $U,V \in \op(M_{sa})$ the sequence
$$\exs{\dbt_M(U \cup
V)}{\dbt_M(U)\oplus\dbt_M(V)}{\dbt_M(U \cap V)}$$ is exact. Then
$\dbt_M$ is a sheaf on $M_{sa}$. Moreover it follows by definition
that $\dbt_M$ is quasi-injective.

\begin{df} Let $Z$ be a closed subset of $M$. We denote by
$\II^\infty_{M,Z}$ the sheaf of $\C^\infty$ functions on $M$
vanishing up to infinite order on $Z$.
\end{df}

\begin{df} A Whitney function on a closed subset $Z$ of $M$ is an indexed
family $F=(F^k)_{k\in \N^n}$ consisting of continuous functions on
$Z$ such that $\forall m \in \N$, $\forall k \in \N^n$, $|k| \leq
m$, $\forall x \in Z$, $\forall \varepsilon >0$ there exists a
neighborhood $U$ of $x$ such that $\forall y,z \in U \cap Z$

$$\left|F^k(z)-\sum_{|j+k|\leq m}{(z-y)^j \over
j!}F^{j+k}(y)\right| \leq \varepsilon d(y,z)^{m-|k|}.$$ We denote
by $W^\infty_{M,Z}$ the space of  Whitney $\C^\infty$ functions on
$Z$. We denote by $\W^\infty_{M,Z}$ the sheaf $U \mapsto
W_{U,U\cap Z}^\infty$.
\end{df}

In \cite{KS96} the authors defined the functor
$$\cdot \wtens \C^\infty_M:\mod_{\rc}(\CC_M) \to \mod(\D_M)$$
in the following way: let $U$ be a subanalytic open subset of $M$
and $Z=M \setminus U$. Then $\CC_U \wtens
\C^\infty_M=\II^\infty_{M,Z}$, and $\CC_Z \wtens
\C^\infty_M=\W^\infty_{M,Z}$. This functor is exact and extends as
a functor in the derived category, from $D^b_{\rc}(\CC_M)$ to
$D^b(\D_M)$. Moreover the sheaf $F \wtens \C^\infty_M$ is soft for
any $\R$-constructible sheaf $F$.

\begin{df} One denotes by $\CWM$ the presheaf of Whitney
$\C^\infty$ functions on $M_{sa}$ defined as follows:
$$U \mapsto \Gamma(M;H^0D'\CC_U \wtens \C^\infty_M).$$
\end{df}
 As a consequence of a result of \cite{Ma67}, for $U,V \in \op(M_{sa})$ the sequence
$$\lexs{\CWM(U \cup
V)}{\CWM(U)\oplus\CWM(V)}{\CWM(U \cap V)}$$ is exact. Then $\CWM$
is a sheaf on $M_{sa}$.

Let us consider a locally cohomologically trivial (l.c.t.)
subanalytic open subset, i.e. $U \in \op({M_{sa}})$ satisfying
$D'\CC_U \simeq \CC_{\overline{U}}$ and $D'\CC_{\overline{U}}
\simeq \CC_U$. Thanks to the triangulation theorem one can prove
that l.c.t. open subanalytic subsets form a basis for the topology
of $M_{sa}$, and given a l.c.t. $U \in \op(M_{sa})$ we have
\begin{eqnarray*}
\Gamma(U;\CWM) & = & \Gamma(M;H^0D'\CC_U \wtens \C^\infty)\\
& \simeq & \Gamma(M;\CC_{\overline{U}} \wtens \C^\infty)\\
& = & W^\infty_{M,\overline{U}}.
\end{eqnarray*}
Moreover $R\Gamma(U;\CWM)$ is concentrated in degree zero since
$\CC_{\overline{U}} \wtens \C^\infty_M$ is soft.\\

Remark that $\Gamma(U;\dbt_M)$ and $\Gamma(U,\CWM)$ are
$\Gamma(\overline{U};\D_M)$-modules for each $U \in \op(M_{sa})$,
hence applying Proposition \ref{6.7.1} the sheaves $\dbt_M$ and $
\CWM$ belong
to $\mod(\rho_!\D_M)$.\\

We have the following result

\begin{prop} For each $F \in D^b_{\rc}(\CC_M)$ one has the isomorphism
\begin{eqnarray*}
\imin \rho\ho(F,\dbt_M) & \simeq & \th(F,\db_M),\\
\imin \rho\rh(F,\CWM) & \simeq & D'F\wtens\C^\infty_M.
\end{eqnarray*}
\end{prop}
\dim\ \ We may reduce to the case $F=k_U$ with $U \in
\op^c(M_{sa})$. Let $V \in \op^c(M_{sa})$.

By definition of $\th$ we have $\Gamma(V;\th(\CC_U,\db_M)) \simeq
\Gamma(U\cap V; \dbt_V)$. Let us consider a subanalytic $W
\subset\subset V$. The natural morphism $\Gamma(U\cap V; \dbt_V)
\to \Gamma(U \cap W;\dbt_M)$ defines the morphism
$$\varphi:\Gamma(U\cap V;\dbt_V) \to \lpro {W\subset\subset V}\Gamma(U\cap W;\dbt_M) \simeq \Gamma(V;\imin \rho \Gamma_U\dbt_M).$$
Since the family $\{W \in \op^c(M_{sa});\; W \subset\subset V\}$
is a covering of $V$ and $\th(\CC_U,\db_M)$ is a sheaf $\varphi$
is an isomorphism.

To prove the second isomorphism we shall first prove the
isomorphism
\begin{equation}\label{CWMwtens}
\ho(F,\CWM) \simeq
H^0D'F\wtens\C^\infty_M
\end{equation}
for $F \in
\mod_{\rc}(\CC_M)$. We may reduce to the case $F=k_U$ with $U$
l.c.t. and subanalytic. Let $V \in \op^c(M_{sa})$ such that $V$
and $U \cap V$ are l.c.t. and let us consider the family $\T=\{W
\in \op^c(M_{sa})\text{ l.c.t.};\; W \subset\subset V,\; W \cap U
\text{ l.c.t.}\}$. The natural morphism $\psi:
\Gamma(V;\CC_{\overline{U}}\wtens\C^\infty_M) \simeq W^\infty_{V,V
\cap \overline{U}} \to W^\infty_{M,\overline{U \cap W}} \simeq
\Gamma(X;\CC_{\overline{U \cap W}}\wtens\C^\infty_M)$ defines the
morphism
$$\psi:\Gamma(V;\CC_{\overline{U}}\wtens\C^\infty_M) \to \lpro {W \in \T}\Gamma(X;\CC_{\overline{U \cap W}}\wtens\C^\infty_M) \simeq \Gamma(V;\imin \rho \Gamma_U\CWM),$$
where the second isomorphism follows since the family $\T$ is
cofinal in $\{W \in \op^c(X_{sa});\; W \subset\subset V\}$. Since
the family $\T$ is a covering of $V$ and
$\CC_{\overline{U}}\wtens\C^\infty_M$ is a sheaf $\psi$ is an
isomorphism. Hence we get the desired isomorphism.

Now let $F \in D^b_{\rc}(\CC_M)$.
We have the chain of morphisms
\begin{eqnarray*}
D'F \wtens \C^\infty_M & \simeq & \indl {F' \to F} \ho(F',\CC_M)
\wtens \C^\infty_M \\
& \simeq & \indl {F' \to F} \imin \rho\ho(F',\CWM) \\
& \to & \indl {F' \to F} \imin \rho\rh(F',\CWM) \\
& \simeq & \imin \rho\rh(F,\CWM),
\end{eqnarray*}
where $F' \to F$ ranges to the family of qis. By Theorem
\ref{eqrhorc} we may suppose $F' \in K^b(\mod_{\rc}(k_X))$ and
then the first isomorphism follows from \eqref{CWMwtens}. We have
$\indl {F' \to F} \ho(F',\CWM) \simeq \rh(F,\CWM)$ and the result
follows since $\imin \rho$ is exact.
\qed \\

Now let $X$ be a complex manifold, $X_\R$ the underlying real
analytic manifold and $\overline{X}$ the complex conjugate
manifold. One denotes by $\ot_X$ and $\OWX$ the sheaves of
tempered and Whitney holomorphic functions respectively which are
defined as follows:
\begin{eqnarray*}
\ot_X:=\rh_{\rho_!\D_X}(\rho_!\OO_{\overline{X}},\dbtxr)\\
\OWX:=\rh_{\rho_!\D_X}(\rho_!\OO_{\overline{X}},\CWXR).
\end{eqnarray*}

By definition, $\ot_X$ and $\OWX$ belong to $D^b(\rho_!\D_X)$. The
relation with the functors of temperate and formal cohomology are
given by the following result 

\begin{prop} For each $F \in D^b_{\rc}(\CC_X)$ one has the
isomorphisms
\begin{eqnarray*}
\imin \rho\rh(F,\ot_X) & \simeq & \th(F,\OO_X),\\
\imin \rho\rh(F,\OWX) & \simeq & D'F \wtens \OO_X.
\end{eqnarray*}
\end{prop}
\dim\ \ We have the chain of isomorphisms
\begin{eqnarray*}
\imin \rho\rh(F,\ot_X) & \simeq & \imin \rho
\rh(F,\rh_{\rho_!\D_X}(\rho_!\OO_{\overline{X}},\dbtxr)) \\
& \simeq & \imin \rho\rh_{\rho_!\D_X}(\rho_!\OO_{\overline{X}},\rh(F,\dbtxr)) \\
& \simeq & \rh_{\D_X}(\OO_{\overline{X}},\imin \rho\rh(F,\dbtxr)) \\
& \simeq & \rh(\OO_{\overline{X}},\th(F,\dbtxr)) \\
& \simeq & \th(F,\OO_X).
\end{eqnarray*}
The proof of $\imin \rho\rh(F,\OWX) \simeq  D'F \wtens \OO_X$ is
similar. \qed


\appendix

\section{Appendix}

\subsection{Review on subanalytic sets}

We recall briefly some properties of subanalytic subsets.
Reference are made to \cite{BM88} and \cite{Lo93}. Let $X$ be a
real analytic manifold.

\begin{df} Let $A$ be a subset of $X$.

\begin{itemize}

\item[(i)] $A$ is said to be semi-analytic if it is locally
analytic, i.e. each $x \in A$ has a neighborhood $U$ such that $X
\cap U=\cup_{i \in I}\cap_{j \in J} X_{ij}$, where $I,J$ are
finite sets and either $X_{ij}=\{y \in U_x;\,f_{ij}>0\}$ or
$X_{ij}=\{y \in U_x;\,f_{ij}=0\}$ for some analytic function
$f_{ij}$.

\item[(ii)] $A$ is said to be subanalytic if it is locally a
projection of a relatively compact semi-analytic subset, i.e. each
$x \in A$ has a neighborhood $U$ such that there exists a real
analytic manifold $Y$ and a relatively compact semi-analytic
subset $A'\subset X \times Y$ satisfying $X \cap U=\pi(A')$, where
$\pi:X \times Y \to X$ denotes the projection.

\item[(iii)] Let $Y$ be a real analytic manifold. A continuous map
$f:X \to Y$ is subanalytic if its graph is subanalytic in $X
\times Y$.
\end{itemize}
\end{df}

Let us recall some result on subanalytic subsets.

\begin{prop} Let $A,B$ be subanalytic subsets of $X$. Then $A \cup
B$, $A \cap B$, $\overline{A}$, $\partial A$ and $A \setminus B$
are subanalytic.
\end{prop}

\begin{prop} Let $A$ be a subanalytic subsets of $X$. Then the
connected components of $A$ are locally finite.
\end{prop}

\begin{prop} Let $f:X \to Y$ be a subanalytic map. Let $A$ be a
relatively compact subanalytic subset of $X$. Then $f(A)$ is
subanalytic.
\end{prop}

\begin{df} A simplicial complex $(K,\Delta)$ is the data
consisting of a set $K$ and a set $\Delta$ of subsets of $K$
satisfying the following axioms:
\begin{itemize}
\item[S1] any $\sigma \in \Delta$ is a finite and non-empty subset
of $K$,
\item[S2] if $\tau$ is a non-empty subset of an element $\sigma$
of $\Delta$, then $\tau$ belongs to $\Delta$,
\item[S3] for any $p \in K$, $\{p\}$ belongs to $\Delta$,
\item[S4] for any $p \in K$, the set $\{\sigma \in \Delta; p \in
\sigma\}$ is finite.
\end{itemize}
\end{df}

If $(K,\Delta)$ is a simplicial complex, an element of $K$ is
called a vertex. Let $\R^K$ be the set of maps from $K$ to $\R$
equipped with the product topology. To $\sigma \in \Delta$ one
associate $|\sigma|\subset \R^K$ as follows:
$$|\sigma|=\left\{x \in \R^K;\text{ $x(p)=0$ for $p \notin
\sigma$, $x(p)>0$ for $p \in \sigma$ and } \sum_p
x(p)=1\right\}.$$ As usual we set:
$$|K|=\bigcup_{\sigma \in \Delta}|\sigma|,$$
$$U(\sigma)=\bigcup_{\tau \in \Delta,\tau \supset \sigma}|\tau|,$$
and for $x \in |K|$:
$$U(x)=U(\sigma(x)),$$
where $\sigma(x)$ is the unique simplex such that $x \in
|\sigma|$.

\begin{teo} Let $X=\bigsqcup_{i \in I} X_i$ be a locally finite
partition of $X$ con\-sisting of subanalytic subsets. Then there
exists a simplicial complex $(K,\Delta)$ and a subana\-lytic
homeomorphism $\psi:|K|\iso X$ such that
\begin{itemize}
\item[(i)] for any $\sigma \in \Delta$, $\psi(|\sigma |)$ is a
subanalytic submanifold of $X$,
\item[(ii)] for any $\sigma \in \Delta$ there exists $i \in I$
such that $\psi(|\sigma |) \subset X_i$.
\end{itemize}
\end{teo}

\subsection{Sheaves on Grothendieck topologies}

We recall the definitions of a Grothendieck topology. We will not
treat the most general case, for which we refer to \cite{KS}. We
will follow the presentation of \cite{Ta94} and \cite{KS01}.

Let $\C$ be a category admitting finite products and fiber
products, and given $U \in \C$, denote by $\C_U$ the category of
arrows $V \to U$. Given a morphism $V \to U$ and $S \subset
\ob(\C_U)$, one denotes by $V \times_U S \subset \ob(\C_V)$ the
subset defined by $\{V \times_U W \to V;\; W \in S\}$.

\begin{df} If, $S_1,S_2 \subset \ob(\C_U)$, one says that $S_1$ is
a refinement of $S_2$ ($S_1 \preceq S_2$ for short) if any $V_1
\to U$ in $S_1$ factorizes as $V_1 \to V_2 \to U$ with $V_2 \to U
\in S_2$.
\end{df}

\begin{df}\label{df:GT} A Grothendieck topology on $\C$ associates to each $U \in \C$ a family $\cov(U) \subset \ob(C_U)$ satisfying the following axioms:

\begin{itemize}

\item[GT1] $\{U \stackrel{id}{\to} U\} \in \cov(U)$,

\item[GT2] if $\cov(U) \ni S_1 \preceq S_2 \subset \ob(\C_U)$, then
$S_2 \in \cov(U)$,

\item[GT3] if $S \in \cov(U)$, then for each $V \to U$, $V \times_U S \in
\cov(V)$,

\item[GT4] if $S_1, S_2 \subset \ob(\C_U)$, $S_1 \in \cov(U)$ and
$V \times_U S_2 \in \cov(V)$, then $S_2 \in \cov(U)$.

\end{itemize}

An object $S \in \cov(U)$ is called a covering of $U$.

\end{df}

\begin{df} A site $X$ is a category $\C_X$ endowed with a Grothendieck topology.

Let $\C_X$ and $\C_Y$ be two categories admitting finite products
and fiber products. A functor of sites $f:X \to Y$ is a functor
$f^t:\C_Y \to C_X$ which commutes with fiber products and such
that if $U \in \C_Y$ and $S \in \cov(U)$, then $f^t(S) \in
\cov(f^t(U))$.

\end{df}

Now let $k$ be a field.

\begin{df} Let $X$ be a site. A presheaf of $k$-modules on $X$ is a functor $\C^{op}_X \to
\mod(k)$.
\end{df}

One denotes by $\psh(k_X)$ the abelian category of presheaves of
$k$-modules on $X$. Let $F \in \psh(k_X)$, let $U \in \C_X$ and
consider $V \to U \in \C_U$. The restriction morphism $F(U) \to
F(V)$ is denoted by $s \mapsto s|_U$.

Let  $F$ be a presheaf of $k$-modules on $X$ and let $S \subset
\ob(\C_U)$. One defines

$$F(S)=\ker\Big(\prod_{V \in S}F(V) \rightrightarrows \prod_{V', V'' \in S}F(V' \times_U V'')\Big)$$

\begin{df} A presheaf $F$ of $k$-modules on $X$
is a separated presheaf (resp. a sheaf) if for each $U \in \C_X$
and each $S \in \cov(U)$ the morphism $F(U) \to F(S)$ is a
monomorphism (resp. an isomorphism).
\end{df}

One denotes by $\mod(k_X)$ the category of sheaves of $k$-modules
on $X$. We set for short $\Ho_{k_X}$ instead of $\Ho_{\mod(k_X)}$.

We recall the construction of a sheaf associated to a presheaf.
The relation ``$\preceq$'' defines a preorder on $\cov(U)$, $U \in
\C_X$. Let $F \in \psh(k_X)$, one defines the functor
$(\cdot)^+:\psh(k_X) \to \psh(k_X)$ in the following way. For each
$U \in \C_X$

$$F^+(U)=\lind{S \in \cov(U)}F(S)$$

\begin{teo}

\begin{itemize}
\item[(i)] The functor $(\cdot)^+:\psh(k_X) \to \psh(k_X)$ is left exact,
\item[(ii)] if $F \in \psh(k_X)$, then $F^+$ is separated,
\item[(iii)] if $F \in \psh(k_X)$ is separated, then $F^+ \in \mod(k_X)$,
\item[(v)] the functor $(\cdot)^{++}:\psh(k_X) \to \mod(k_X)$ is
exact,
\item[(iv)] let $F \in \psh(k_X)$ and $G \in \mod(k_X)$, one has the adjunction formula:
$$\Ho_{\psh(k_X)}(F,\iota G) \simeq \Ho_{k_X}(F^{++},G),$$
where $\iota$ denotes the embedding functor.
\end{itemize}

\end{teo}

Let $F \in \psh(k_X)$, the sheaf $F^{++}$ is called the sheaf
associated to $F$.

\begin{prop} Let $F,G \in \mod(k_X)$. A morphism $\varphi \in
\Ho_{k_X}(F,G)$ is an epimorphism if and only if for each $U \in
\C_X$ there exists $\{U_i\}_{i \in I} \in \cov(U)$ such that for
each $s \in G(U)$ there exists $t_i \in F(U_i)$ such that
$\varphi(t_i)=s|_{U_i}$ for each $i$.
\end{prop}

Let $f:X \to Y$ be a morphism of sites. Let $F \in \psh(k_X)$ and
$G \in \psh(k_Y)$. One defines the functors

\begin{eqnarray}
f_*: \psh(k_X) & \to & \psh(k_Y) \label{eq:f*} \\
f^{\gets}:\psh(k_Y) & \to & \psh(k_X)
\end{eqnarray}
in the following way: let $U \in \C_X$ and $V \in \C_Y$, then

\begin{eqnarray*}
f_*F(V) & = & F(f^t(V)) \\
f^{\gets}F(U) & = & \lind {U \to f^t(W)}G(W),
\end{eqnarray*}
where $W \in \C_Y$.

\begin{df} Let $f:X \to Y$ be a functor of sites

\begin{itemize}
\item[(i)] the functor of direct image $f_*: \mod(k_X) \to \mod(k_Y)$ is the functor induced by
(\ref{eq:f*}),
\item[(ii)] the functor of inverse image $f^{-1}: \mod(k_Y) \to \mod(k_X)$ is defined by
$f^{-1}=(f^{\gets}(\cdot))^{++}$.
\end{itemize}

\end{df}

\begin{prop}
\begin{itemize}
\item[(i)] The functor $f_*$ is left exact and commutes with
$\Lpro$,
\item[(ii)] the functor $f^{-1}$ is exact and commutes with
$\Lind$,
\item[(iii)] $(\imin f, f_*)$ is a pair of adjoint functors.
\end{itemize}

\end{prop}

\begin{df} Let $X$ be a site and let $F,G \in \mod(k_X)$.

\begin{itemize}
\item[(i)] One denotes by $\ho(F,G)$ the sheaf $U \mapsto
\Ho_{k_U}(F|_U,G|_U)$,
\item[(ii)] one denotes by $F \otimes G$ the sheaf associated to the presheaf $U \mapsto F(U) \otimes
G(U)$.

\end{itemize}

\end{df}

\begin{prop} Let $F \in \mod(k_X)$, $G,G' \in \mod(k_Y)$.

\begin{itemize}

\item[(i)] $\ho_{k_Y}(G,f_*F) \simeq f_*\ho_{k_X}(f^{-1}G,F)$,
\item[(ii)] $f^{-1}(G \otimes G') \simeq f^{-1}G \otimes
f^{-1}G'$.
\end{itemize}

\end{prop}

\addcontentsline{toc}{section}{

\end{document}